\documentclass[10pt]{elsarticle}
\makeatletter
\def\ps@pprintTitle{%
 \let\@oddhead\@empty
 \let\@evenhead\@empty
 \def\@oddfoot{\centerline{\thepage}}%
 \let\@evenfoot\@oddfoot}
\makeatother
\usepackage{amsmath,amssymb,amsthm,amsfonts,amsopn,natbib,mathrsfs}

\usepackage{color,geometry,graphicx,dsfont}
\usepackage{url,setspace}
\usepackage[usenames,dvipsnames,svgnames,table]{xcolor}
\usepackage{caption,subcaption}

\theoremstyle{plain}
\newtheorem{theorem}{Theorem}[section]
\newtheorem{proposition}[theorem]{Proposition}
\newtheorem{lemma}[theorem]{Lemma}
\newtheorem{corollary}[theorem]{Corollary}
\theoremstyle{definition}
\newtheorem{definition}[theorem]{Definition}
\theoremstyle{remark}
\newtheorem{remark}[theorem]{Remark}
\newtheorem{example}[theorem]{Example}

\newcommand{\bd}{\mathbf} 
\newcommand{\dsi}{\mathds{1}} 

\newcommand{\LtR}{\mathbf L^2(\mathbb{R})}
\newcommand{\LiR}{\mathbf L^1(\mathbb{R})}
\def\esssup{\mathop{\operatorname{ess~sup}}}
\def\essinf{\mathop{\operatorname{ess~inf}}}
\newcommand{\abs}[1]{|#1|}
\newcommand{\RR}{\mathbb R}
\newcommand{\ZZ}{\mathbb Z}
\newcommand{\NN}{\mathbb N}
\newcommand{\CC}{\mathbb C}
\newcommand{\TT}{\mathbb T}
\newcommand{\inner}[2]{\langle #1,#2 \rangle}
\newcommand{\norm}[1]{\|#1\|}
\newcommand{\sgn}{\mathop{\operatorname{sgn}}}  
\newcommand{\supp}{\mathop{\operatorname{supp}}} 
\newcommand{\revised}[1]{{#1}}

\DeclareFontFamily{U}{mathx}{\hyphenchar\font45}
\DeclareFontShape{U}{mathx}{m}{n}{
      <5> <6> <7> <8> <9> <10>
      <10.95> <12> <14.4> <17.28> <20.74> <24.88>
      mathx10
      }{}
\DeclareSymbolFont{mathx}{U}{mathx}{m}{n}
\DeclareFontSubstitution{U}{mathx}{m}{n}
\DeclareMathAccent{\widecheck}{0}{mathx}{"71}
\DeclareMathAccent{\wideparen}{0}{mathx}{"75}
\newcommand{\LD}{\mathbf L^2(D)}
\newcommand{\LtDF}{\mathcal F^{-1}(\mathbf L^2(D))}

\begin{document}
 
 
 \begin{frontmatter}
 \title
 {A class of warped filter bank frames tailored to non-linear frequency scales}
 
 \author[ari]{Nicki Holighaus\corref{cor1}}
 \ead{nicki.holighaus@oeaw.ac.at}
 
 \author[ait]{Christoph Wiesmeyr}
 \ead{christoph.wiesmeyr@ait.ac.at}
 
 \author[ari]{Zden\v{e}k Pr\r{u}\v{s}a}
 \ead{zdenek.prusa@oeaw.ac.at}
 
 \cortext[cor1]{Corresponding author}
 
 \address[ari]{Acoustics Research Institute Austrian Academy of Sciences, Wohllebengasse 12--14, A-1040 Vienna, Austria}
 \address[ait]{AIT Austrian Institute of Technology GmbH, Donau-City-Strasse 1, A-1220 Vienna, Austria and\\ NuHAG, Faculty of Mathematics, University of Vienna, Oskar-Morgenstern-Platz 1, A-1090 Vienna, Austria}

  \begin{abstract}
   A method for constructing non-uniform filter banks is presented.
   Starting from a uniform system of translates, generated by a prototype filter, a non-uniform covering of the frequency axis is obtained by composition with a warping function. The warping function is a $\mathcal C^1$-diffeomorphism that determines the frequency progression and can be chosen freely, apart from minor technical restrictions. The resulting functions are interpreted as filter frequency responses. Combined with appropriately chosen decimation 
   factors, a non-uniform analysis filter bank is obtained. Classical Gabor and wavelet filter banks are special cases of the proposed construction.
   Beyond the state-of-the-art, we construct a filter bank adapted to a frequency scale derived from human auditory perception and families of filter banks that can be interpreted as 
   an interpolation between linear (Gabor) and logarithmic (wavelet) frequency scales. We derive straightforward conditions on the prototype filter decay and the decimation factors, such that the resulting warped filter bank forms a frame.    
   In particular, a simple and constructive
   method for obtaining tight frames with bandlimited filters is derived by invoking previous results on generalized shift-invariant systems. 
 \end{abstract}
 

\end{frontmatter}

 \section{Introduction}
   
   In this contribution, we introduce a class of non-uniform time-frequency systems optimally adapted to non-linear frequency scales. The central paradigm of our construction, and what distinguishes it from previous approaches, is to provide uniform frequency resolution \emph{on the target frequency scale}. Invertible time-frequency systems are of particular importance, since they allow for stable recovery of signals from the time-frequency representation coefficients. Therefore, we also    
   derive necessary and sufficient conditions for the resulting systems to form a frame. 
   
   To demonstrate the flexibility and importance of our construction, illustrative examples recreating (or imitating) classical time-frequency representations such as Gabor~\cite{ga46,gr01,fest98,fest03}, wavelet~\cite{ma09-1,da92} or $\alpha$-transforms~\cite{cofe78,fo89,hona03} are provided. While this paper considers the setting of (discrete) Hilbert space frames, the properties of continuous warped time-frequency systems are investigated in the related contribution~\cite{bahowi14}
   . Whenever a time-frequency filter bank adapted to a given frequency progression and with linear time-progression in each channel is desired, we believe that the proposed \emph{warped filter banks} provide the right framework for its design.
   
   Generalized shift-invariant (GSI)  systems~\cite{rosh05,helawe02,chhale15,jale14,badohojave11} over $\LtR$ are families $(g_{m,n})_{m,n\in\ZZ}$ with $g_{m,n} = \widetilde{g_m}(\cdot - na_m)$, for filters  $(\widetilde{g_m})_{m\in\ZZ}\subset\LtR$ and decimation factors $(a_m)_{m\in\ZZ}\subset\RR^+$. In the proposed method, GSI systems are constructed from a prototype frequency response $\theta$ via composition with a \emph{warping function} $\Phi$ that specifies the desired frequency scale/progression. \revised{Applying the inverse Fourier transform $\mathcal{F}^{-1}$, we obtain $\widetilde{g_m} = \mathcal{F}^{-1}(\theta(\cdot - m) \circ \Phi)$} up to normalization. To highlight the relation of the resulting warped time-frequency systems to non-uniform filter banks, we use terminology from filter bank theory and refer to GSI systems as filter banks.
   
   \textbf{Contribution.} We connect the frame theory of abstract GSI systems with warped filter banks. In particular, we use the structure of warped filter banks to show the following:
   \begin{itemize}
    \item[(i)] Under mild restrictions, warped filter banks satisfy the important \emph{local integrability condition} for GSI systems, see \cite{helawe02,chhale15}.
    \item[(ii)] The application of results from \cite{ho14-1,chhale15} provides intuitive necessary Bessel and frame conditions for warped filter banks that are easy to verify.
    \item[(iii)] If $\theta$ is compactly supported and the decimation factors $(a_m)_{m\in\ZZ}\subset\RR^+$ are small enough, these necessary conditions are also sufficient, we are in the \emph{painless case}~\cite{badohojave11}. In this setting, the canonical dual of a warped filter bank frame is a warped filter bank as well.
    \item[(iv)] If $\theta$ has sufficient decay, but not necessarily compact support, then there exist decimation factors $(a_m)_{m\in\ZZ}\subset\RR^+$ that yield warped filter bank frames. It is discussed how the notion of sufficient decay depends on the warping function $\Phi$. 
   \end{itemize}
   Most of the above results are made possible by choosing $(a_m)_{m\in\ZZ}\subset\RR^+$ majorized by a set of \emph{natural decimation factors}. Natural decimation factors are obtained by observing that the bandwidth (or essential Fourier domain support) of the $g_{m,n}$ is intrinsically linked to the derivative of the warping function $\Phi$. For natural decimation factors, the selection of a single parameter $\tilde{a}>0$ determines all decimation factors $a_m$ in a way that respects the bandwidth of the warped filters. We further provide examples for \emph{tight warped filter banks} from compactly supported $\theta$ and a construction of discrete warped filter banks for digital signals in $\ell^2(\ZZ)$. The latter is adapted from~\cite{hoprwi15}, where warped filter banks for $\ell^2(\ZZ)$ were first presented. For discrete warped filter banks we experimentally verify that, the frame bound ratio deteriorates slowly when the decimation requirements of the painless case are violated. Our results \revised{indicate} that a frame bound ratio below $10$ can be achieved with very little oversampling.
   
   \textbf{Adapted time-frequency systems.} Time-frequency (or time-scale) representations are an indispensable tool for signal analysis and processing. The most widely used and most thoroughly explored such representations  are certainly Gabor and wavelet transforms and their variations, e.g. windowed modified cosine~\cite{brpr86,brjopr87} or wavelet packet~\cite{comequwi94,Wickerhauser:2027800} transforms. The aforementioned transforms  unite two very important properties: There are various, well-known necessary and/or sufficient conditions for stable inversion from the transform coefficients, i.e., for the generating function system to form a frame. In addition to the perfect reconstruction property, the frame property ensures stability of the synthesis operation after coefficient modification, enabling controlled time-frequency processing. Furthermore, efficient algorithms for the computation of the transform coefficients and the synthesis operation exist for each of the mentioned transforms~\cite{st98-1,ma09-1}.
 
   While providing a sound and well-understood mathematical foundation, Gabor and wavelet transforms are designed to follow two specified frequency scales: linear, respectively logarithmic.  A wealth of approaches exists to soften this restriction, e.g. decompositions using filter banks \cite{Cvetkovic:2003a,Cvetkovic:1998a,Vaidyanathan:1993a,bofehl98}, for example based on perceptive frequency scales \cite{ha97-X,Strahl:2009a,Patterson:1992a}. Adaptation over time is considered in approaches such as modulated lapped transforms~\cite{ma92-2}, adapted local trigonometric transforms~\cite{wewi93-1} or (time-varying) wavelet packets~\cite{heorraxi97}. Techniques that jointly offer flexible time-frequency resolution and variable redundancy, the perfect reconstruction property and efficient computation are scarce however. The setting of so-called nonstationary Gabor transforms~\cite{badohojave11}, a recent generalization of classical Gabor transforms, provides the latter $2$ properties while allowing for freely 
chosen time  progression and varying resolution. In this construction, the 
   frequency scale is still linear, but the sampling density may be changed over time. The properties of nonstationary Gabor systems have been investigated in, e.g., \cite{ho14-1,doma14-1,doma14-2}. When desiring increased flexibility along frequency, generalized shift-invariant systems~\cite{rosh05,helawe02,chel05-X,ch16-X,chhale15}, or equivalently (non-uniform) \emph{filter banks}~\cite{akva03}, provide the analogous concept. They offer full flexibility in frequency, with a linear time progression in each filter, but flexible sampling density across the filters. Analogous, continuously indexed systems are considered in \cite{basp14,jale14}. Indeed, nonstationary Gabor systems are equivalent to filter banks via an application of the (inverse) Fourier transform to the generating functions. Note that all the widely used transforms mentioned in the  previous paragraph can be interpreted as filter banks.
 
   \textbf{Adaptation to non-linear frequency scales through warping.}
   \revised{There have been previous attempts to construct adapted filter banks by frequency warping. In our approach, a warping operator is simply a coordinate change and not unitary in general. In contrast, all previous methods consider unitary warping following the change of variables formula for integration. Thus, they introduce an additional weight function depending on $\Phi$ and modify the filter values beyond a simple coordinate change. Consequently, the properties of both types of warping, the resulting systems and the challenges faced in their construction are quite different.} 
   
   For example, Braccini and Oppenheim~\cite{brop74}, as well as Twaroch and Hlawatsch~\cite{hltw98}, propose a unitary warping of a \revised{system} of translates, interpreted as filter frequency responses. In~\cite{brop74} only spectral analysis is desired, while time-frequency distributions are constructed in~\cite{hltw98}, without considering signal reconstruction.
   
   The application of unitary warping to an entire Gabor or wavelet system has also been investigated~\cite{bajo93-1,ba94-2,caev98,doevma13}. Although unitary transformation bequeaths basis (or frame) properties to the warped atoms, the resulting system is not anymore a filter bank. Instead, the warped system produces undesirable, dispersive time-shifts and the resulting representation is not easily interpreted, see \cite{caev98}. Only for the continuous short-time Fourier transform, or under quite strict assumptions on a Gabor system, a \emph{redressing} procedure can be applied to recover a GSI system~\cite{evangelista2013warped}. In all other cases, the combination of unitary warping with redressing complicates the efficient, exact computation of redressed warped Gabor frames, such that approximate implementations are considered~\cite{evangelista2014approximations}.
 
   Finally, it should be noted that the idea of a (non-unitary) logarithmic warping of the frequency axis to obtain wavelet systems from a system of translates was already used in the proof of the so called \emph{painless conditions} for wavelets systems~\cite{dagrme86}. Recent parallel work by Christensen and Goh \cite{chgo14} focuses on exposing the duality between Gabor and wavelet systems via the mentioned logarithmic warping. However, the idea has never been relaxed to other frequency scales so far. In the present work, we generate time-frequency transformations beyond wavelet and Gabor systems by allowing more general warping functions. The proposed warping procedure has already proven useful in the area of graph signal processing~\cite{shuman2015spectrum}. 
 


 \section{Preliminaries}\label{sec:pre}

We use the following normalization of the Fourier transform
  $\hat f(\xi) := \mathcal{F}f(\xi) = \int_\RR f(t) e^{-2 \pi i t \xi} \; dt$,
  for all $f\in\LiR$
and its unitary extension to $\LtR$. The inverse Fourier transform is denoted by $\widecheck{f} = \mathcal F^{-1}f$. For an open interval $D\subset \RR$, typically $D = \RR$ or $D=\RR^+$, we use the convention that $\mathbf{L}^2(D) := \{f \in\mathbf{L}^2(\RR) \colon  f(t) = 0 \text{ for almost every } t\in\RR\setminus D\}$, such that the Fourier transform and its inverse restrict naturally to $\mathbf{L}^2(D)$.
Following this convention, we denote by 
\[ \LtDF  \subseteq \LtR\]
the space of functions whose Fourier spectrum is restricted to $D$.

Further, we frequently use the \emph{translation operator} defined by $\bd T_x f = f(\cdot -x)$, for all $f \in \LtR$, 
and the composition $f\circ g := f(g(\cdot))$ of two functions $f$ and $g$. The standard Lebesgue measure is denoted by $\mu$.

When discussing the properties of the constructed function systems in the following sections, we will repeatedly use the notions of weight functions and weighted $\mathbf{L}^p$-spaces, $1\leq p\leq \infty$. Weighted $\mathbf{L}^p$-spaces are defined as
\begin{equation*}
  \mathbf{L}^p_w(\RR) := \left\{ f:\RR\mapsto \CC \colon  wf\in \mathbf{L}^p(\RR) \right\}.
\end{equation*}
with a continuous, \revised{positive function $w:\RR\mapsto \RR^+$ called \emph{weight function}. The associated norm is $\|f\|_{\mathbf{L}^p_w}:= \|wf\|_{\mathbf{L}^p}$. In the following, when the term weight function is used, continuity and positivity are always implied.}

Two special classes of weight functions are of particular interest: Let $v: \RR \rightarrow \RR^+$ and $w: \RR \rightarrow \RR^+$ be \emph{continuous, positive weight functions}. We call $v: \RR \rightarrow \RR^+$ \emph{submultiplicative} and $w: \RR \rightarrow \RR^+$ \emph{$v$-moderate} respectively if they satisfy, for some $C>0$, 
\begin{equation}\label{eq:moderateness}
	v(x+y)\leq v(x)v(y), \text{ and } w(x+y) \leq C v(x) w(y),\quad \text{ for all } x,y\in\RR.
\end{equation}
In particular, we can (and will) always choose $v$ such that $1$ is a valid choice for the constant in the latter inequality ($\max\{C,1\}v$ is submultiplicative  whenever $v$ is).
Submultiplicative and moderate weight functions play an important role in the theory of function spaces, as they are closely related to the translation-invariance of the corresponding weighted spaces \cite{fe79-2,gr01}, see also~\cite{gr07} for an in-depth analysis of weight functions and their role in harmonic analysis.

A generalized shift-invariant (GSI) system on $\LtR$ is a union of shift-invariant systems $\{\bd T_{na_m} h_m \in\LtR  \colon  n\in\ZZ \}$, with $h_m\in\LtR$ and $a_m\in\RR^+$, for all $m$ in some index set. The representation coefficients of a function $f\in\LtR$ with respect to the GSI system are given by the inner products
\[
 c_f(n,m):= \langle f, \bd T_{na_m} h_m\rangle = \left(f\ast \overline{h_m(-\cdot)}\right)(na_m), 
\]
for all $n,m$. Here, we denote by $h_m(-\cdot)$ the map $t\mapsto h_m(-t)$. The above representation of the coefficients in terms of a convolution \revised{of $f$ with the conjugate, time-inverse of $h_m$} alludes to the fact that $c_f(\cdot,m)$ is a filtered, and sampled, version of $f$. This relation justifies our use of filter bank terminology when discussing GSI systems. 

\begin{definition}\label{def:filterbank}
  Let $\left( g_{m} \right)_{m\in \ZZ}\subset \LD$  and $\left( a_m \right)_{m\in \ZZ} \subset \RR^+$. 
  We call the system
  \begin{equation}
    \left( g_{m,n} \right)_{m,n\in\ZZ},\quad g_{m,n}:=\bd T_{na_m} \mathcal{F}^{-1}(g_m), \text{ for all } n,m\in\ZZ,
    \label{}
  \end{equation}
  a \emph{(non-uniform) filter bank} for $\LtDF$.
  The elements of $\left( g_{m} \right)_{m\in \ZZ}$  and $\left( a_m \right)_{m\in \ZZ}$ are called \emph{frequency responses} and \emph{decimation factors}, respectively.
\end{definition}

Such filter banks can be used to analyze signals in $\LtDF$, and for a given signal $f\in \LtDF$, we refer to the sequence $c_f := (c_f(n,m))_{n,m\in\ZZ} = \left( \inner{f}{g_{m,n}} \right)_{m,n \in \ZZ} $ as the \emph{filter bank (analysis) coefficients}. A \emph{uniform filter bank} is a filter bank with $a_m = a$ for all $m\in\ZZ$. 

For many applications it is of great importance that all the considered signals can be reconstructed from these coefficients, in a stable fashion. It is a central observation of frame theory that this is equivalent to the existence of constants $0<A\leq B<\infty $, such that
\begin{equation}\label{eq:frameprop}
  A \norm{f}_2^2 \leq \norm{c_f}_{\ell^2(\ZZ^2)}^2 \leq B \norm{f}_2^2, \; \text{for all $f\in \LtDF$}.
\end{equation}
A system $\left( g_{m,n} \right)_{m,n\in\ZZ}$ that satisfies this condition is called \emph{filter bank frame}~\cite{dusc52,ch16-X}. A (filter bank) frame \revised{is called \emph{tight} if equality can be achieved in \eqref{eq:frameprop} and \emph{snug} if $B/A\approx 1$ is possible}. If at least the upper inequality in \eqref{eq:frameprop} is satisfied, then $\left( g_{m,n} \right)_{m,n\in\ZZ}$ is a \emph{Bessel sequence}. In that case, the frame operator is defined by 
\begin{equation}\label{eq:frameop}
  \bd S: \LtDF \rightarrow \LtDF, \quad\bd{S} f = \sum_{m,n\in\ZZ} c_f(m,n) g_{m,n}, \text{ for all } f\in\LtDF.
\end{equation}
If $\left( g_{m,n} \right)_{m,n\in\ZZ}$ is a frame, then the frame operator is invertible. It is the key component in the construction of the \emph{canonical dual frame} $\left( \widetilde{g_{m,n}} \right)_{m,n\in \ZZ}$, obtained by applying the inverse of the frame operator to the frame elements, i.e., $\widetilde{g_{m,n}} := \bd S^{-1}(g_{m,n})$, for all $m,n \in \ZZ$. The canonical dual frame facilitates \emph{perfect reconstruction} from the analysis coefficients:
\begin{equation}
  f = \sum_{m,n \in \ZZ} c_f(m,n) \widetilde{g_{m,n}}, \text{ for all } f\in \LtDF. 
\end{equation}
Note that, in contrast to short-time Fourier or uniform filter bank frames, there is no guarantee that the canonical dual frame, or indeed any dual frame, of a general filter bank frame is of the form $( \bd T_{na_m} \mathcal{F}^{-1}(\widetilde{g_m}) )_{n,m\in\ZZ}$, for some $(\widetilde{g_m})_{m\in\ZZ}\subset \LD$ and the same sequence of decimation factors $\left( a_m \right)_{m\in \ZZ}$.
Abstract filter bank frames~\cite{bofehl98} have received considerable attention, as (generalized) shift-invariant systems in \cite{ja98,helawe02,rosh05,chhale15,jale14} and as (frequency-side) nonstationary Gabor systems in \cite{badohojave11,doma14-1,doma14-2,ho14-1}. In contrast, this contribution is concerned with a specific, structured family of filter bank systems and how the superimposed structure can be used to construct filter bank frames. 
    

\section{Warped filter banks}\label{sec:warp}

  In signal analysis, the usage of different frequency scales has a long history. Linear and logarithmic scales arise naturally when constructing a filter bank through modulation or dilation of a single prototype filter, respectively. In this way, the classical Gabor and wavelet transforms are obtained. The consideration of alternative frequency scales can be motivated, for example, from (a) theoretical interest in a family of time-frequency representations that serve as an interpolation between the two extremes, as is the case for the $\alpha$-transform\revised{~\cite{cofe78,fo89,hona03}} (which can be related to polynomial scales), or (b) specific applications and/or signal classes. A prime example for the second case is audio signal processing with respect to an auditory frequency scale, e.g. in gammatone filter banks~\cite{Patterson:1992a,Strahl:2009a} adapted to the ERB scale~\cite{glmo90}, the latter modeling the frequency progression and frequency-bandwidth relationship in the human cochlea. The mentioned 
  methods have several things in common: They are all based on a single prototype filter and possess the structure of a GSI system (or filter bank). The bandwidth of their filters is directly linked to the filter center frequencies and their spacing, induced by the frequency scale. 

    The filter banks we propose in this section have the property that they are designed as a system of translates on a given frequency scale. 
    This scale determines a conversion from frequency to a new unit (e.g. ERB) with respect to which the designed filters provide a uniform resolution.
In the next sections, we will show that this construction admits a special class of non-uniform filter banks with a simplified structure compared to general filter banks.
    
  Formally, a frequency scale is specified by a continuous, bijective function $\Phi: D\rightarrow \RR$ and the transition between the non-linear scale $\Phi$ and the linear scale is achieved by $\Phi$ and $\Phi^{-1}$. Hence, we construct filter frequency responses from a prototype function $\theta:\RR\mapsto\CC$ by means of \emph{translation}, followed by \emph{deformation}, 
  \begin{equation}\label{eq:thetafm}
    \left( (\bd T_m\theta)\circ \Phi \right)_{m\in \ZZ}.
  \end{equation}
  This general formulation provides tremendous flexibility for frequency scale design. Furthermore,  choosing $\Phi$ as $\Phi(\xi)\mapsto a\xi$ or $\Phi(\xi)\mapsto \log_a(\xi)$, for $a>0$, yields systems of translates $\bd T_{m/a}\left(\theta(a\cdot)\right)$ and dilates $(\theta\circ \log_a)(\cdot/a^m)$, respectively. 
  Such $\Phi$ will provide the starting point for recovering Gabor and wavelet filter banks in our framework. 
 
\begin{definition}\label{def:warpfun}
    Let $D\subseteq \RR$ be any open interval. A $\mathcal C^1$-diffeomorphism $\Phi: D \rightarrow \RR$ is called \emph{warping function}, if 
    \begin{enumerate}
     \item[(i)] the derivative $\Phi'$ of $\Phi$ is positive, i.e., $\Phi'>0$, and
     \item[(ii)] there is a submultiplicative weight $v$, such that the weight function 
     \begin{equation}\label{eq:weightFunDef}
       w: = \left( \Phi^{-1} \right)' = \frac{1}{\Phi'\left(\Phi^{-1}(\cdot) \right)}
     \end{equation}
     is $v$-moderate, i.e., $w(\tau_0+\tau_1)\leq v(\tau_0)w(\tau_1)$, for all $\tau_0,\tau_1\in\RR$.     
    \end{enumerate}
  \end{definition}
  
  Given a warping function $\Phi$, $w$ and $v$ will from now on always denote weights as specified in Definition \ref{def:warpfun}. 
  
  \begin{remark}
    While moderateness and invertibility of $\Phi$ will prove essential for our results, there are no technical obstructions preventing us from allowing warping functions $\Phi\in\mathcal C^0(D)\setminus \mathcal C^1(D)$, such that $\Phi'$ is only piecewise continuous. However, this implies that some (or all) of the elements of the warped family given in  \eqref{eq:thetafm} can have at most piecewise continuous derivative, independent of the smoothness of $\theta$ and with the implied negative effects to their Fourier localization. Moreover, $\Phi$ is easily lifted to $\mathcal C^1$ with minor, arbitrarily local changes. Therefore, we only see limited value in generalizing the notion of a warping function beyond diffeomorphisms.
  \end{remark}
  
  \begin{proposition}\label{pro:warpfun_dil_mul}
    If $\Phi: D \rightarrow \RR$ is a warping function as per Definition \ref{def:warpfun}, then $\tilde{\Phi}:=c\Phi(\cdot/d)$ is a warping function with domain $dD$, for all positive, finite constants $c,d\in\RR^+$. If $w = (\Phi^{-1})'$ is $v$-moderate, then $\tilde{w}=(\tilde{\Phi}^{-1})'$ is $v(\cdot/c)$-moderate.
  \end{proposition}
  \begin{proof}
      The result is easily obtained by elementary manipulation.
  \end{proof}

  Several things should be noted when considering the definition and proposition above.
  \begin{itemize}
   \item Proposition \ref{pro:warpfun_dil_mul} shows that it really is sufficient to consider integer translates of the prototype $\theta$ when constructing the frequency responses $\theta_{\Phi,m}$. If $a>0$ is arbitrary, then with $\theta_a := \theta(\cdot/a)$, we have 
\begin{equation}
(\bd T_m \theta_a)\circ (a\Phi) = \theta_a\left(a\Phi(\cdot)-m\right) = \theta\left(\Phi(\cdot)-m/a\right) = (\bd T_{m/a}\theta)\circ\Phi.
\label{eq:scadil}
\end{equation}
   \item Moderateness of $w = \left( \Phi^{-1} \right)'$ ensures translation invariance of the associated weighted $\mathbf{L}^p$-spaces. In particular, identifying $(\bd T_m\theta)\circ \Phi$ with its trivial extension to the whole real line, we have
  \begin{equation}\label{eq:moderatenormbound}
    \norm{(\bd T_m\theta)\circ \Phi}^2_{\LD} = \norm{\bd T_m \theta}^2_{\mathbf L^2_{\sqrt{w}}(\RR)} \leq \begin{cases}
      v(m) \norm{\theta}^2_{\mathbf L^2_{\sqrt{w}}(\RR)} &,\ \text{if } \theta \in \mathbf L^2_{\sqrt{w}}(\RR)\\                                                                                           
      w(m) \norm{\theta}^2_{\mathbf L^2_{\sqrt{v}}(\RR)} &,\ \text{if } \theta \in \mathbf L^2_{\sqrt{v}}(\RR).                                                                                            
    \end{cases}    
  \end{equation}
  \item $\mathbf L^2_{\sqrt{v}}(\RR)\subseteq \mathbf L^2_{\sqrt{w}}(\RR)$, since \eqref{eq:moderatenormbound}, with $m=0$, implies $\norm{\theta}^2_{\mathbf L^2_{\sqrt{w}}(\RR)}\leq w(0)\norm{\theta}^2_{\mathbf L^2_{\sqrt{v}}(\RR)}$.
 \end{itemize}

  A warped filter bank can now be constructed easily. To do so, after selecting the warping function $\Phi$, one simply chooses an appropriate prototype frequency response $\theta$ and positive decimation factors $(a_m)_{m\in\ZZ}$.
  
  \begin{definition}\label{def:warpedfilterbanks}
  Let $\Phi: D\rightarrow \RR$ be a warping function and $\theta \in \mathbf L^2_{\sqrt{v}}(\RR)$. Furthermore, let $\mathbf a := \left( a_m \right)_{m\in \ZZ} \subset \RR^+$ be a set of decimation factors. Then the \emph{warped filter bank} with respect to the triple $(\Phi,\theta,\mathbf a)$ is given by
  \begin{equation}\label{eq:warpedfilterbank}
    \mathcal G(\Phi,\theta,\mathbf a) := \left( \bd T_{n a_m} \widecheck{g_m} \right)_{m,n \in \ZZ} = \left( \bd T_{n a_m} \mathcal F^{-1}(g_m) \right)_{m,n \in \ZZ},
  \end{equation}
  with
  \begin{equation}
    g_m(\xi) := \begin{cases} 
             \sqrt{a_m} (\bd T_m \theta)\circ \Phi(\xi) & \text{ if } \xi\in D,\\
             0 & \text{ else.}
           \end{cases}
    \label{eq:gmDefintion}
  \end{equation}
  If $a_m = \tilde a/w(m)$, for all $m\in\ZZ$ and some $\tilde{a}>0$, then we say that $\mathbf a$ is a set of \emph{natural decimation factors (for $\Phi$)}.
 \end{definition}
 
  Although, in theory, the choice of decimation factors is arbitrary, it is worth emphasizing the importance \revised{of} natural decimation factors: If a set of decimation factors $\mathbf a$ is majorized by some natural decimation factors, then the frequency responses $g_m$ are uniformly $\bd L^2$-bounded, recall \eqref{eq:moderatenormbound}. \revised{Uniform} $\bd L^2$-boundedness is easily seen to be a necessary condition for the Bessel property, i.e., the existence of an upper bound in \eqref{eq:frameprop}. Moreover, we will see in Section \ref{sec:warpedframes} that natural decimation factors play an essential role in the creation of warped filter bank frames.
  
  In Gabor theory, the set of lattices that provide certain properties, e.g., the frame property, is often studied for a fixed window. Similarly, for warped filter banks, we might fix the warping function $\Phi$ and the prototype filter $\theta$. Then, we can attempt to determine decimation factors $\mathbf a$, such that $\mathcal G(\Phi,\theta,\mathbf a)$ is a frame or possess some other property of interest. In that sense, the decimation factors $\mathbf a$ should be interpreted as parameters of the warped filter bank that are chosen and \emph{tuned} post-hoc to achieve the desired properties.   

\begin{remark}
 Note that the condition $\theta \in \mathbf L^2_{\sqrt{w}}(\RR)$ is sufficient to ensure that $\mathcal G(\Phi,\theta,\mathbf a) \subset \LtDF$. In that setting, a set of natural decimation factors would 
 have the form $a_m = \tilde a/v(m)$ for all $m\in\ZZ$ and some $\tilde{a}>0$, instead. All results in this contribution also hold in this case and are proven with the same techniques. Since a decay condition on $\theta$ is usually considered less severe than a restriction of the sampling density, our results are presented for the configuration given in Definition \ref{def:warpedfilterbanks}. 
 
 However, depending on how much the submultiplicative weight $v$ deviates from $w=(\Phi^{-1})'$, the two sets of natural decimation factors, and the spaces of eligible prototype functions, may differ significantly. Therefore, we shortly discuss the necessary changes in the case $\theta \notin \mathbf L^2_{\sqrt{v}}(\RR)$ in Section \ref{sec:ThetaInLtw}.
\end{remark}

Assume for now that $\theta\in \mathbf L^2_{\sqrt{v}}(\RR) \cap \mathbf L^1_{w}(\RR)$. If we rewrite the elements of $\mathcal G(\Phi,\theta,\mathbf a)$ as a Fourier integral, i.e., 
\[
 \begin{split}
 \mathcal F^{-1}(g_m)(t) & = \int_{\RR} g_m(\xi) e^{2\pi i \xi t}~d\xi\\
 & = \sqrt{a_m} \int_{D} \theta(\Phi(\xi)-m) e^{2\pi i \xi t}~d\xi = \sqrt{a_m} \int_{\RR} w(\tau+m)\theta(\tau) e^{2\pi i \Phi^{-1}(\tau+m) t}~d\tau,
 \end{split}
\]
with the change of variable $\xi = \Phi^{-1}(\tau + m)$, we can see that decay (smoothness) of $\theta$ implies smoothness (decay) for the elements of $\mathcal G(\Phi,\theta,\mathbf a)$, provided that $\Phi$ is smooth enough as well. This behavior is crucial for the construction of systems with good time-frequency localization and in fact central for the results presented in~\cite{bahowi14}, where the above Fourier integrals are studied in more detail.

We now provide some examples of warping functions that are of particular interest, e.g., because they encompass important frequency scales. In Proposition \ref{pro:warpingfunctions} at the end of this section, we show that the presented examples indeed define warping functions in the sense of Definition \ref{def:warpfun}. Some instances of the warping functions in the following examples can be seen in Figure \ref{fig:warpingfuncts}.

{
\begin{example}[Wavelets]\label{ex:wavelet1}
  Choosing $\Phi=\log$, with $D=\RR^+$ leads to a system of the form
  \begin{equation*}
    g_{m}(\xi)= \sqrt{a_m}\theta(\log(\xi)-m)=\sqrt{a_m}\theta(\log(\xi e^{-m})) = \sqrt{\frac{a_m}{a_0}}g_0(\xi e^{-m}).
    \label{}
  \end{equation*}
  This warping function therefore leads to $g_{m}$ being a dilated version of some $g_{0} = \theta\circ\log$. The natural decimation factors are given by $a_m = \tilde a/w(m)= \tilde a e^{-m}$.
  This shows that $\mathcal G(\log, \theta, \tilde a e^{-m})$ is indeed a wavelet system, with the minor modification that our scales are reciprocal to the usual definition of wavelets.
\end{example}

\begin{example}\label{ex:Rplusalpha1}
  The family of warping functions $\Phi_l(\xi) = c\left((\xi/d)^l - (\xi/d)^{-l}\right)$, for some $c,d >0$ and $l \in (0,1]$, is an alternative to the logarithmic warping for the domain $D=\RR^+$. The logarithmic warping in the previous example can be interpreted as the limit of this family for $l\rightarrow 0$ in the sense that for any fixed $\xi\in \RR^+$,
  \begin{equation}\label{eq:limitArg}
    \Phi_l'(\xi)/l = \frac{c}{d}\left((\xi/d)^{-1+l} + (\xi/d)^{-1-l}\right) \overset{l\rightarrow 0}{\rightarrow} \frac{2c}{\xi} = \frac{2c}{d}\log'(\xi/d).
  \end{equation}

  This type of warping provides a frequency scale that approaches the limits $0$ and $\infty$ of the frequency range $D$ in a slower fashion than the wavelet warping. In other words, $g_{m}$ is less deformed for $m>0$, but more deformed for $m<0$ than in the case $\Phi=\log$. 
  Furthermore, the property that $g_{m}$ can be expressed as a dilated version of $g_{0}$ is lost.  
\end{example}

\begin{example}[ERBlets]\label{ex:ERB1}
  In psychoacoustics, the investigation of filter banks adapted to the spectral resolution of the human ear has been subject to a wealth of research, see \cite{momo03} for an overview.
  We mention here the Equivalent Rectangular Bandwidth scale (ERB-scale) described in \cite{glmo90}, which introduces a set of bandpass filters modeling human perception, see also \cite{bahoneso13} for the construction of an invertible filter bank adapted to the ERB-scale. In our terminology the ERB warping function is given by
  \begin{equation*}
      \Phi_{\text {ERB}}(\xi) = \sgn{(\xi)} \; c \log\left(1+\frac{\abs{\xi}}{d}\right), \\
    \label{}
  \end{equation*}
  where the constants are given by $c=9.265$ and $d=228.8$. \revised{Using this function, we obtain a filter bank with uniform bandwidth \emph{on the ERB-scale}. This property is the key feature of the filter banks proposed in \cite{bahoneso13} or more traditional Gammatone filter banks on the ERB-scale, see e.g. \cite{Strahl:2009a}. However, in contrast to these previous methods, it is very easy to construct tight warped filter bank frames, see Section \ref{sec:warpedframes}.} The warped ERB filter bank has potential applications in audio signal processing, as it provides an invertible transform adapted to the human perception of sound.
\end{example}

\begin{example}\label{ex:alpha1}
  Filter banks obtained from the warping functions $\Phi_\alpha(\xi) = \sgn(\xi)  \,\left((\abs{\xi}+1)^{1-\alpha}-1\right)$, for some $\alpha \in [0,1)$ can serve as a substitute for the \emph{$\alpha$-transform}, see \cite{cofe78,fo89,hona03,fefo06}. The latter is a filter bank constructed from a single prototype frequency response by translation and dilation, leading to \revised{filters with bandwidth proportional to $(1+|\xi|)^\alpha$}, where $\xi\in\RR$ is the center frequency of the filter frequency response. Varying $\alpha$, one can \emph{interpolate} between the Gabor transform ($\alpha = 0$, constant \revised{bandwidth}) and a wavelet-like (or more precisely ERB-like) transform with the \revised{bandwidth} depending linearly on the center frequency ($\alpha\rightarrow 1$). It is easy to confirm that the warping function $\Phi_\alpha(\xi)$ yields \revised{filters with bandwidth proportional to $(1+|\xi|)^\alpha$ as well}. We will see in subsequent sections that, in stark contrast to the $\alpha$-transform, it is easy to 
construct tight frames using the warping function $\Phi_\alpha$. Note as well that study of the $\alpha$-transform usually excludes the limiting case $\alpha=1$, which is also not captured by the above warping construction. However, the logarithmic warping considered in Example \ref{ex:ERB1} yields filters with bandwidth proportional to $(1+|\xi|)$ and can thus be considered as substitute for the limiting case. \revised{A limit argument on the derivative of $\Phi_\alpha$ ($\alpha\rightarrow 1$), similar to \eqref{eq:limitArg}, provides another approach to confirm this assertion.}
\end{example}

\begin{example}\label{ex:tan}
  Finally, we propose a warping function for representing functions band-limited to the interval $D=(-\pi,\pi)$. For this purpose set $\Phi(\xi) = \tan(\xi)$. Necessarily, the frequency responses $g_m$, given by \eqref{eq:gmDefintion}, are all compactly supported on $D$ and increasingly peaky and concentrated at the upper and lower borders of $D$, as $m$ tends to $\infty$ and $-\infty$, respectively. By using the equivalence of GSI systems and nonstationary Gabor systems~\cite{badohojave11} through application of the Fourier transform, we can thus construct time-frequency systems on arbitrary open intervals. Frames for intervals have been proposed previously by Abreu et al.~\cite{Abreu2014274}.
\end{example}

\begin{figure}[!tbh]
  \includegraphics[width=.5\textwidth]{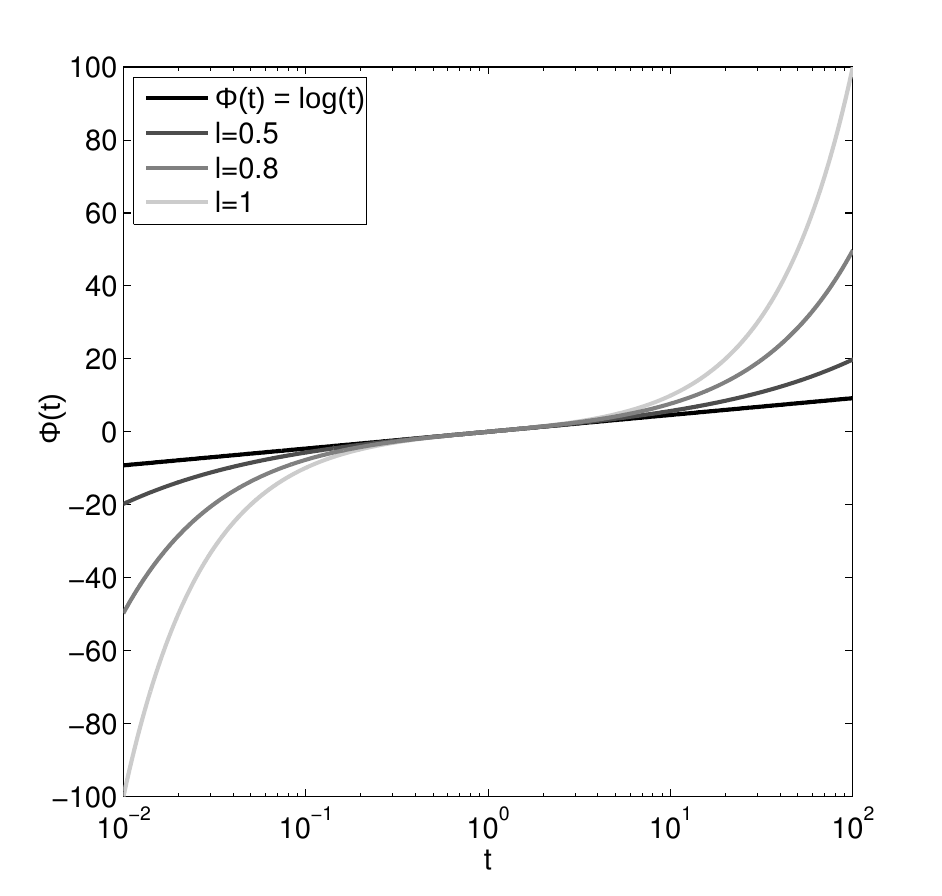}\includegraphics[width=.5\textwidth]{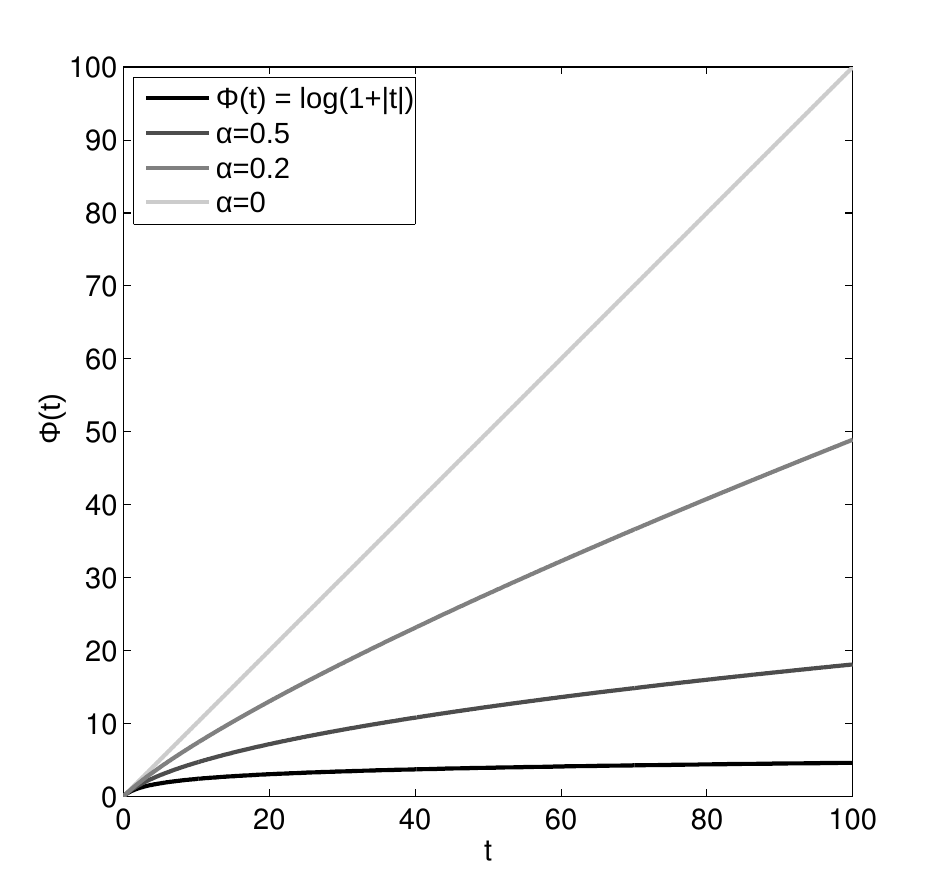}
  \caption{(left) Warping functions from Examples \ref{ex:wavelet1} and \ref{ex:Rplusalpha1}: This plot shows logarithmic (wavelet) warping function (black) and $\Phi_l = l^{-1}(\xi^l-\xi^{-l})$, for $l=0.5$ (dark gray), $l=0.8$ (medium gray) and $l=1$ (light gray). Note that the horizontal axis is logarithmic. (right) Warping functions from Examples \ref{ex:ERB1} and \ref{ex:alpha1}. This plot shows the ERBlet warping function (black), with $c=d=1$,  and $\Phi_\alpha = (1-\alpha)^{-1}\sgn(\xi)((1+|\xi|)^{(1-\alpha)}-1)$, for $\alpha=0.5$ (dark gray), $\alpha=0.2$ (medium gray) and $\alpha=0$ (light gray). The horizontal axis is linear.}\label{fig:warpingfuncts}
\end{figure}

\begin{proposition}\label{pro:warpingfunctions}
  The following are valid triples of warping functions, weights $w$ and moderating submultiplicative weights $v$, as per Definition \ref{def:warpfun}:
  \begin{enumerate}
    \item[(i)] $\Phi: \RR^+ \rightarrow \RR, \; \xi \mapsto \log(\xi)$, with $w = v = e^{(\cdot)}$.
    \item[(ii)] $\Phi: \RR^+ \rightarrow \RR, \; \xi \mapsto \left((\xi)^{l} - (\xi)^{1-l}\right)$, for $l \in (0,1]$, with 
    \[ 
    w = \frac{2^{-1/l}}{l}\cdot\frac{((\cdot)+\sqrt{(\cdot)^2+4})^{1/l}}{\sqrt{(\cdot)^2+4}}\ \text{ and }\ v = \sqrt{2(\cdot)^2+8}\cdot\left(\sqrt{\frac{5}{4}}\cdot \left(2+\frac{|\cdot|}{2}\right)\right)^{2/l}.
    \]
    \item[(iii)] $\Phi: \RR \rightarrow \RR, \; \xi \mapsto \sgn(\xi) \, \log(1+\abs{\xi})$ with $w = v = e^{|\cdot|}$.
    \item[(iv)] $\Phi: \RR \rightarrow \RR, \; \xi \mapsto \sgn(\xi)  \, \left((1+\abs{\xi})^{1-\alpha}-1\right)$, for some $\alpha \in [0,1)$, with $w = v/(1-\alpha) = (1-\alpha)^{-1}(1+|\cdot|)^{\alpha/(1-\alpha)}$.
    \item[(v)] $\Phi: (-\pi,\pi) \rightarrow \RR, \; \xi \mapsto \tan(\xi)$, with $w = (1+(\cdot)^2)^{-1}$ and $v = 2(1+(\cdot)^2)$.
  \end{enumerate}  
\end{proposition}

\begin{proof}
  Items (i), (iii) and (iv) are easily shown through elementary calculations. It remains to prove items (ii) and (v).
  
  Ad (ii): $\Phi$ is in $\mathcal{C}^\infty(\mathbb{R}^+)$ with $\Phi^{-1}(\tau_0) = 2^{-1/l}(\tau_0+\sqrt{\tau_0^2+4})^{1/l}$, such that 
  \begin{equation*}
    w(\tau_0) = \left( \Phi^{-1} \right)'(\tau_0) = \frac{2^{-1/l}}{l}\cdot\frac{(\tau_0+\sqrt{\tau_0^2+4})^{1/l}}{\sqrt{\tau_0^2+4}} = \frac{\Phi^{-1}(\tau_0)}{l}(\tau_0^2+4)^{-1/2}.
  \end{equation*}  
  Thus, $w$ is $v$-moderate with $v= v_0v_1$, if $\Phi^{-1}$ is $v_0$-moderate and $((\cdot)^2+4)^{-1/2}$ is $v_1$-moderate. It is straightforward to show that the latter is true with $v_1(\tau_0) = \sqrt{2}\sqrt{\tau_0^2+4}$. We now proceed to show that $\Phi^{-1}$ is $v_0$-moderate with $v_0(\tau_0) = (5/4)^{1/l}(2+|\tau_0|/2)^{2/l}$. Observe that 
  \begin{equation*}
   \left(\frac{\Phi^{-1}(\tau_0+\tau_1)}{\Phi^{-1}(\tau_0)}\right)^{l}-1 =    \frac{\tau_1+\sqrt{(\tau_1+\tau_0)^2+4}-\sqrt{\tau_0^2+4}}{\tau_0+\sqrt{\tau_0^2+4}} \leq \begin{cases} 
                                                                                                   0 & \text{ for all } \tau_1\leq 0, \\
                                                                                                   \frac{\tau_1}{\sqrt{\tau_1^2/4+4}-\tau_1/2} & \text{ else,}
                                                                                                \end{cases}
  \end{equation*}
  since the left hand side attains its global maximum at $\tau_0=-\tau_1/2$, for fixed $\tau_1>0$. By the 
  fundamental theorem of calculus (with $f = \sqrt{\tau_1^2/4+(\cdot)}$), 
  \begin{equation*}
    \frac{\tau_1}{\sqrt{\tau_1^2/4+4}-\tau_1/2} \leq \frac{\tau_1^2+16}{4} \leq (2+|\tau_1|/2)^2, 
  \end{equation*}
  valid for all $\tau_1>0$. Hence, 
  \begin{equation*}
     \left(\frac{\Phi^{-1}(\tau_0+\tau_1)}{\Phi^{-1}(\tau_0)}\right)^{l} \leq 1 + (2+|\tau_1|/2)^2 \leq 5/4 \cdot (2+|\tau_1|/2)^2
  \end{equation*}
  and the right hand side is a submultiplicative weight. Consequently, $\Phi^{-1}$ is indeed $v_0$-moderate and we obtain that $w$ is $v$-moderate with 
  \begin{equation*}
    v(\tau_0) = v_0(\tau_0)v_1(\tau_0) = \sqrt{2\tau_0^2+8}\cdot\left(\sqrt{\frac{5}{4}}\cdot \left(2+\frac{|\tau_0|}{2}\right)\right)^{2/l},
  \end{equation*}
  as desired.

  Ad (v): The crucial step is to show that $\arctan' = (1+(\cdot)^2)^{-1}$ can be moderated by a submultiplicative weight. But since $(\tau_0 + \tau_1)^2 \leq 2(\tau_0^2+\tau_1^2)$ for all $\tau_0,\tau_1\in\RR$, it is easy to see that $v = 2(1+(\cdot)^2)$ is submultiplicative and that $w=(1+(\cdot)^2)^{-1}$ is $v$-moderate. The other required properties of $\Phi = \tan$ are elementary.
\end{proof}
}


 \section{Warped filter bank frames}\label{sec:warpedframes}
 
 Although this contribution is concerned only with warped filter bank frames, our results are derived from structural properties and results obtained in the general, abstract filter bank (or GSI) setting. As such, the structure imposed on warped filter banks can be seen as a constructive means to satisfy, or simplify, the conditions of these abstract results. Our results rely on the simple, but crucial identity
 \begin{equation}\label{eq:warpedfreqresp}
    \sum_{m\in\ZZ} a_m^{-1}|g_m(\xi)|^2 = \sum_{m\in\ZZ} |(\bd T_m\theta)\circ \Phi(\xi)|^2, \text{ for all } \xi\in D,
 \end{equation}
 a direct consequence of the definition of the warped filter bank $\mathcal G(\Phi,\theta,\bd a)$. As a consequence of the above equality, we can find upper and lower bounds for  \eqref{eq:warpedfreqresp}, by instead determining upper and lower bounds on the simpler quantity
 \begin{equation}\label{eq:freqrespWarp}
   \sum_{m\in\ZZ} |\bd T_m \theta|^2.
 \end{equation}
 
 \revised{To exclude} pathological cases from the study of filter bank frames, it has proven useful to assume that a filter bank $(g_{m,n})_{m,n\in\ZZ}$ satisfies the so-called local integrability condition~\cite{helawe02,jale14,chhale15}. This enables the generalization of numerous important results, e.g. a characterization of dual frames, from the frame theory of Gabor systems~\cite{gr01} and uniform filter banks~\cite{ja98}. 
  
 \begin{definition}\label{def:LICandALIC}
   Denote by $\mathcal D$ the set of all functions $f\in\bd L^\infty(D)$ with compact support. We say that the filter bank $(g_{m,n})_{m,n\in\ZZ}$, generated from 
   $\left( g_{m} \right)_{m\in \ZZ}\subset \LD$  and $\left( a_m \right)_{m\in \ZZ} \subset \RR^+$, satisfies the local integrability condition (LIC), if 
   \begin{equation}\label{eq:LIC}
     L(f):= \sum_{m\in\ZZ}\sum_{l\in\ZZ} a_m^{-1} \int_{\supp(f)} \left|f(\xi+la_m^{-1})g_m(\xi)\right|^2~d\xi <\infty,
   \end{equation}
   for all $f\in \mathcal D$.
 \end{definition}
 
 The LIC might seem intimidating and opaque at first, but once we impose some structure on $(g_{m,n})_{m,n\in\ZZ}$, it can often be substituted by mild conditions on the frequency responses $g_m$ and decimation factors $a_m$. In the case of warped filter banks, boundedness of \eqref{eq:freqrespWarp} and $\bd a$ being majorized by a set of natural decimation factors is already sufficient for $\mathcal G(\Phi,\theta,\mathbf a)$ to satisfy the LIC.
 
 \begin{theorem}\label{thm:warpedLIC}
  Let $\Phi:D\rightarrow \RR$ be a warping function and $\theta\in \bd L^2_{\sqrt{v}}(\RR)$. If 
  \begin{equation}\label{eq:suffLICcond}
    \sup_{m\in\ZZ} a_m w(m)<\infty\qquad \text{and}\qquad \esssup_{\tau\in\RR} \sum_{m\in\ZZ} |\bd T_m \theta(\tau)|^2 < \infty,
  \end{equation}
  then $\mathcal G(\Phi,\theta,\mathbf a)$ satisfies the LIC \eqref{eq:LIC}. In particular, if $\bd a$ is a set of natural decimation factors and the second condition in \eqref{eq:suffLICcond} holds, then $\mathcal G(\Phi,\theta,\mathbf a)$ satisfies the LIC.
 \end{theorem}
 
 \begin{proof}
  First note that, instead of considering all compactly supported and essentially bounded functions $f\in\LD$, it is sufficient to verify the LIC \eqref{eq:LIC} only for the characteristic functions $\dsi_I$ on compact intervals $I\subset D$. These functions are clearly contained in $\LD$ and it is easy to see that 
  \[
   \supp(f) \subseteq I\quad \Longrightarrow\quad L(f) \leq \|f\|_{\infty}^2 L(\dsi_I).
  \]
  For $\dsi_I$, the LIC reads
  \begin{equation}\label{eq:charLIC}
    L(\dsi_I) = \sum_{m\in\ZZ} a_m^{-1} \sum_{l\in\ZZ} \int_I \dsi_{I+la_m^{-1}}(\xi) \left|g_m(\xi)\right|^2~d\xi. 
  \end{equation}
  If the right hand side of \eqref{eq:charLIC} is finite, then it converges absolutely and we can interchange sums and integrals freely. Hence, 
  \[
    \begin{split}
      L(\dsi_I) & = \sum_{m\in\ZZ} a_m^{-1} \int_I \left|g_m(\xi)\right|^2 \sum_{l\in\ZZ} \dsi_{I+la_m^{-1}}(\xi) ~d\xi\\
      & < \sum_{m\in\ZZ} \frac{a_m\mu(I) + 1}{a_m} \int_I \left|g_m(\xi)\right|^2~d\xi\\
      & = \sum_{m\in\ZZ} (a_m\mu(I) + 1) \int_I \left|(\bd T_m \theta)\circ \Phi(\xi)\right|^2~d\xi,
    \end{split}
   \]
   where we used that $\sum_{l\in\ZZ} \dsi_{I+la_m^{-1}}(\xi) \leq \lceil a_m\mu(I)\rceil < a_m\mu(I)+1$ for arbitrary $\xi\in D$. We split the upper estimate into two terms and interchange integration and summation once more to obtain
   \[
    L(\dsi_I) < \int_I \sum_{m\in\ZZ} \left|(\bd T_m \theta)\circ \Phi(\xi)\right|^2~d\xi + \sum_{m\in\ZZ} a_m\mu(I)\cdot \int_I \left|(\bd T_m \theta)\circ \Phi(\xi)\right|^2~d\xi.
   \]
   By assumption, there is some constant $B>0$, such that $\sum_{m\in\ZZ} |\bd T_m \theta(\tau)|^2 < B$ almost everywhere and we can conclude that 
   \[
    \int_I \sum_{m\in\ZZ} \left|(\bd T_m \theta)\circ \Phi(\xi)\right|^2~d\xi \leq \mu(I)B.
   \]
   To estimate the second term, note that the change of variable $\xi = \Phi^{-1}(\tau+m)$ yields
   \[     
       \sum_{m\in\ZZ} a_m \cdot \int_I \left|(\bd T_m \theta)\circ \Phi(\xi)\right|^2~d\xi = \sum_{m\in\ZZ} a_m \cdot \int_{\Phi(I)-m} w(\tau+m) \left|\theta(\tau)\right|^2~d\tau = (\ast).
   \]
   By assumption $\sup_{m\in\ZZ} a_m w(m)<\infty$ and $\theta\in \bd L^2_{\sqrt{v}}(\RR)$. To estimate the right hand side of the above equation, we can 
   use $v$-moderateness of $w$:
   \[
    \begin{split}
      (\ast) & \leq \sum_{m\in\ZZ} a_m w(m) \cdot \int_{\Phi(I)-m} v(\tau)\left|\theta(\tau)\right|^2~d\tau\\
      & = \sum_{m\in\ZZ} a_m w(m) \cdot \int_\RR \dsi_{\Phi(I)-m} v(\tau)\left|\theta(\tau)\right|^2~d\tau\\
      & = \sup_{m\in\ZZ} a_m w(m) \cdot \int_\RR v(\tau)\left|\theta(\tau)\right|^2 \sum_{m\in\ZZ} \dsi_{\Phi(I)-m}~d\tau\\
      & < \left(1+\mu(\Phi(I))\right) \cdot  \sup_{m\in\ZZ} a_m w(m) \cdot \int_\RR v(\tau)\left|\theta(\tau)\right|^2 ~d\tau\\
      & = \left(1+\mu(\Phi(I))\right)\cdot \sup_{m\in\ZZ} a_m w(m) \cdot \|\theta\|^2_{\bd L^2_{\sqrt{v}}} < \infty.
    \end{split}
   \]
   Altogether, we obtain 
   \[
    L(\dsi_I) < \mu(I)\cdot\left( B + \left(1+\mu(\Phi(I))\right)\cdot \sup_{m\in\ZZ} a_m w(m) \cdot \|\theta\|^2_{\bd L^2_{\sqrt{v}}}\right) < \infty,
   \]
   which establishes the desired result.
   If $\bd a$ is a set of natural decimation factors, then $a_m w(m) = \tilde a < \infty$ for all $m\in\ZZ$, yielding the second claim. 
\end{proof}
 
 \revised{The previous result shows that a large class of warped filter banks satisfies the LIC. In that case, the results presented in \cite{helawe02,jale14,chhale15}, many of which require the LIC, are fully applicable. Moreover, we will see in Theorem \ref{thm:PLwarp} that the second condition in \eqref{eq:suffLICcond} is in fact necessary for $\mathcal G(\Phi,\theta,\mathbf a)$ to be a Bessel sequence. In other words, for a warped filter bank with the Bessel property, the LIC can always be satisfied by choosing appropriate decimation factors. 
 Our next result relies on some previous results from the literature, which we now recall. Their application will yield necessary and sufficient conditions for a warped filter bank to form a Bessel sequence or even a frame.}
 
 \begin{proposition}\label{pro:necandsuf}   
  Let $(g_{m,n})_{m,n\in\ZZ}$ the filter bank generated from $(g_m)_{m\in\ZZ}\subset \LD$ and $(a_m)_{m\in\ZZ}\subset \RR^+$.
  \begin{itemize}
   \item[(i)] \cite[Prop. 3]{ho14-1} If $(g_{m,n})_{m,n\in\ZZ}$ is a Bessel sequence with bound  $B<\infty$, then 
    \begin{equation}\label{eq:boundeddiag}
	\sum_{m\in\ZZ} \frac{1}{a_m} |g_m(\xi)|^2 \leq B, \text{ for almost all $\xi \in D$. }
    \end{equation}
    \item[(ii)] \cite[Cor. 3.4]{chhale15} If $(g_{m,n})_{m,n\in\ZZ}$ is a frame with lower frame bound $A>0$ satisfying the LIC \eqref{eq:LIC}, then 
     \begin{equation}\label{eq:boundeddiagbel}
       A \leq \sum_{m\in\ZZ} \frac{1}{a_m} |g_m(\xi)|^2, \text{ for almost all $\xi \in D$.}
     \end{equation}     
    \item[(iii)] \cite[Cor. 1]{badohojave11} Assume that there are some constants $c_m,d_m\in\RR$, such that $\supp(g_m) \subseteq [c_m,d_m]$ and
      $a_m$ satisfies $a_m^{-1} \geq d_m-c_m$, for all $m\in\ZZ$. Then $(g_{m,n})_{m,n\in\ZZ}$ forms a frame, with frame bounds $A,B$, for $\LtDF$ if and only if 
    \begin{equation}\label{eq:boundeddiagpl}
	0 < A \leq \sum_{m\in\ZZ} \frac{1}{a_m} |g_m(\xi)|^2 \leq B < \infty, \text{ for almost all } \xi \in D.
    \end{equation}
    Furthermore, the filter bank generated from $(\widetilde{g_m})_{m\in\ZZ}\subset \LD$ and $(a_m)_{m\in\ZZ}\subset \RR^+$, with
    \begin{equation}\label{eq:pldual}
      \widetilde{g_m} = \frac{g_m}{\sum_{l\in\ZZ} \frac{1}{a_l} |g_l|^2}, \text{ for all } m\in\ZZ,
    \end{equation}
    is the canonical dual frame for $(g_{m,n})_{m,n\in\ZZ}$.
  \end{itemize}
\end{proposition}

 With the above results in place, we obtain the following necessary and sufficient conditions for warped filter bank frames.
 
 \begin{theorem}\label{thm:PLwarp}   
  Let $\mathcal{G}(\Phi,\theta,\bd{a})$ be a warped filter bank for $\LtDF$.
  \begin{itemize}
   \item[(i)] If $\mathcal{G}(\Phi,\theta,\bd{a})$ is a Bessel sequence with bound $B<\infty$, then 
    \begin{equation}
	\sum_{m\in\ZZ}  |\bd T_m \theta(\tau)|^2 \leq B < \infty, \text{ for almost all } \tau \in \RR.\label{eq:boundeddiag2a}
    \end{equation}   
   \item[(ii)] If $\mathcal{G}(\Phi,\theta,\bd{a})$ is a frame with lower bound $A>0$ and $\sup_{m\in\ZZ} a_m w(m)<\infty$, then 
    \begin{equation}
	0 < A \leq \sum_{m\in\ZZ}  |\bd T_m \theta(\tau)|^2, \text{ for almost all } \tau \in \RR.
	\label{eq:boundeddiag2b}
    \end{equation}
  \item[(iii)] Assume that there are constants $c<d$, such that $\supp(\theta) \subseteq [c,d]$ and $a_m^{-1} \geq \Phi^{-1}(d+m)-\Phi^{-1}(c+m)$, for all $m\in\ZZ$. The warped filter bank $\mathcal{G}(\Phi,\theta,\bd{a})$ forms a frame for $\LtDF$, with frame bounds $A,B$, if and only if $0< A\leq \sum_{m\in\ZZ}  |\bd T_m \theta|^2 \leq B < \infty$ almost everywhere.
  Furthermore, the canonical dual frame for $\mathcal{G}(\Phi,\theta,\bd{a})$ is given by $\mathcal{G}(\Phi,\tilde \theta,\bd{a})$, with
  \begin{equation}\label{eq:pldual2}
      \tilde{\theta} = \frac{\theta}{\sum_{l\in\ZZ} |\bd T_l \theta|^2}.
  \end{equation}
  \end{itemize}
\end{theorem}
\begin{proof}
   Part (i) is a direct consequence of Proposition \ref{pro:necandsuf}(i) and \eqref{eq:warpedfreqresp}. Part (ii) follows similarly from Proposition \ref{pro:necandsuf}(ii) and \eqref{eq:warpedfreqresp}, after noting that the assumptions of (ii) imply (i) and thus $\sup_{m\in\ZZ} a_m w(m)<\infty$ yields the LIC by Theorem \ref{thm:warpedLIC}. Finally, (iii) is obtained by inserting \eqref{eq:warpedfreqresp} into Proposition \ref{pro:necandsuf}(iii).
\end{proof}

Note that the canonical dual frame in Theorem \ref{thm:PLwarp}(iii) is a warped filter bank as well, obtained with the warping function $\Phi$ and decimation factors $\bd a$. \revised{The dual prototype filter $\tilde \theta$ is easily computed using Equation \eqref{eq:pldual2}.} Theorem \ref{thm:PLwarp}(iii) is the natural generalization of the classical \emph{painless nonorthogonal expansions}~\cite{dagrme86} to warped filter banks and extremely useful when strictly bandlimited filters are required. The whole of Theorem \ref{thm:PLwarp} serves as a strong indicator that for any \emph{snug} frame, i.e., with $B/A\approx 1$, the sum $\sum_{m\in\ZZ}  |\bd T_m \theta|^2$ must necessarily be close to constant. It is thus imperative that the translates of the original window $\theta$ have good summation properties. \revised{We established an intimate relationship between stability of the filter bank $\mathcal{G}(\Phi,\theta,\bd{a})$ and the extrema of $\sum_{m\in\ZZ}  |\bd T_m \theta|^2$.}

Sometimes, when we would like to work in the setting of Theorem \ref{thm:PLwarp}(iii) it can be more efficient to estimate the support of the $g_m$ instead of calculating it exactly. The following result and its discussion below show that this can easily be done using natural decimation factors that satisfy the conditions of Theorem \ref{thm:PLwarp}(iii). The derived natural decimation factors \revised{often} satisfy $a_m^{-1} \approx \Phi^{-1}(d+m)-\Phi^{-1}(c+m)$.

For the purpose of the following result and for later use, we define the function $V:\RR\times \RR \rightarrow \RR$ by
\begin{equation}\label{eq:defOfV}
 V(\tau_0,\tau_1) := \int_{\tau_0}^{\tau_1} v(\tau)~d\tau, \text{ for all } \tau_0,\tau_1\in\RR.
\end{equation}

\begin{corollary}\label{cor:PLnatural}
  Let $\mathcal G(\Phi,\theta,\bd a)$ be a warped filter bank with compactly supported prototype $\theta\in \bd L^2(\RR)$. Define 
  \[
      c_0 := \inf\, \supp(\theta)\quad \text{and}\quad d_0 := \sup\, \supp(\theta)
  \]
  and $\widetilde{a_w}:= V(c_0,d_0)^{-1}$. Assume furthermore that 
  \[
   a_m \leq \widetilde{a_w}/w(m),\quad \text{for all } m\in\ZZ.
  \]
  The warped filter bank $\mathcal G(\Phi,\theta,\bd a)$ forms a frame with frame bounds $A,B$, if and only if $0< A\leq \sum_{m\in\ZZ}  |\bd T_m \theta|^2 \leq B < \infty$ almost everywhere. In that case, the canonical dual frame is given by $\mathcal G(\Phi,\tilde \theta,\bd{a})$, with $\tilde{\theta}$ as in  \eqref{eq:pldual2}.
\end{corollary}

Note that $V(c_0,d_0)$ can be bounded from above by $(d_0-c_0)\max_{\tau\in[c_0,d_0]} v(\tau)$. This coarser estimate can be used for an even simpler computation of decimation factors appropriate for Corollary \ref{cor:PLnatural}, e.g. if $v$ is nondecreasing away from zero.

\begin{proof}[Proof of Corollary \ref{cor:PLnatural}]
  By  the fundamental theorem of calculus, 
  \[
   \begin{split}
   \Phi^{-1}(d_0+m)-\Phi^{-1}(c_0+m) &= \int_{c_0+m}^{d_0+m} w(\tau)~d\tau\\
    & \leq w(m)\cdot \int_{c_0}^{d_0} v(\tau)~d\tau = V(c_0,d_0) w(m),
    \end{split}
  \]
  for all $m\in\ZZ$. Therefore, 
  \[
   \Phi^{-1}(d_0+m)-\Phi^{-1}(c_0+m) \leq w(m)/\widetilde{a_w} \leq a_m^{-1},
  \]
  as per the assumption. Since $\theta\in\LtR$ with compact support implies $\theta\in\bd L^2_{\sqrt{v}}(\RR)$, we can apply Theorem \ref{thm:PLwarp}(iii) to finish the proof.
\end{proof}

Without additional assumptions on the warping function $\Phi$, the condition $a_m \leq \widetilde{a_w}/w(m)$ in Corollary \ref{cor:PLnatural} cannot be improved.
To see this we construct a warping function, such that any choice $a_m>\widetilde{a_w}/w(m)$ yields $a_m^{-1}<\Phi^{-1}(d)-\Phi^{-1}(c)$, for all nonempty, closed intervals $[c,d]$. Hence, the conditions of Theorem \ref{thm:PLwarp}(iii) are violated.
Choose $\Phi = \log$ and note $\Phi^{-1}(\tau) = e^\tau = w(\tau)$ for all $\tau\in\RR$. We can choose $v = w$ and obtain 
\[
 e^{d_0+m}-e^{c_0+m} = e^m\int_{c_0}^{d_0} e^\tau d\tau = V(c_0,d_0)w(m) = (\Phi^{-1}(d_0)-\Phi^{-1}(c_0))v(m),
\]
to show that for a logarithmic warping function the natural decimation factors are indeed the coarsest possible decimation factors that satisfy the conditions Theorem \ref{thm:PLwarp}(iii).

 \subsection{On tight warped filter bank frames}\label{sec:tight}
 
 In the following, we will demonstrate how the definition of warped filter banks leads to straightforward constructions of tight frames with compactly supported prototype $\theta$, following the decimation conditions in Theorem \ref{thm:PLwarp}(iii). \revised{Specifically, the tight frame property is achieved by} selecting a prototype $\theta$ that is compactly supported and equalizes the inequality \eqref{eq:boundeddiag2a}. Tight frames are important for various reasons. They provide a perfect reconstruction system in which \revised{analysis and synthesis frame are equal up to a constant}. Hence, there is no need for computing and/or storing a dual frame, which might be highly inefficient. By equalizing the frame inequality \eqref{eq:frameprop}, they provide optimal norm stability. Furthermore, the usage of tight frames guarantees that the synthesis shares the properties of the analysis, e.g. in terms of time-frequency localization. 

Under the conditions of Theorem \ref{thm:PLwarp}(iii), i.e., $\theta$ is compactly supported and the $a_m$ are small enough, a warped filter bank $\mathcal G(\Phi,\theta,\bd a)$ is a tight frame if and only if, for some $C>0$, 
  \begin{equation}\label{eq:POUsqrt}
    \sum_{m\in \ZZ} \abs{\bd T_m \theta}^2 = C, \text{ a.e.}
  \end{equation}  
However, even if the conditions of Theorem \ref{thm:PLwarp}(iii) are not satisfied, \eqref{eq:POUsqrt} is still a necessary condition for the frame property, at least if the decimation factors $\bd a$ are majorized by a set of natural decimation factors. Therefore, $\theta$ that satisfy  \eqref{eq:POUsqrt} are the optimal starting point when aiming \revised{to construct snug warped filter bank frames, i.e., frames with small frame bound ratio.}

Although surely not the only methods for obtaining functions satisfying  \eqref{eq:POUsqrt}, we highlight here two classical methods that provide both compact support, which is required to apply Theorem \ref{thm:PLwarp}(iii), and a prescribed smoothness: B-splines~\cite{de78} and windows constructed as a superposition of truncated cosine waves of different frequency~\cite{nu81-1}. The second class contains classical window functions such as the Hann, Hamming and Blackman windows.
We now recall a procedure to construct such functions that also satisfy  \eqref{eq:POUsqrt}. The method has previously been reported and proven as \cite[Theorem 1]{shuman2015spectrum}:

  Let $K\in {\mathbb N}$ and $c_k \in \RR$ for $k \in \{0,1,\ldots,K\}$, and define
  \begin{equation}\label{eq:window_form}
   \vartheta(\tau) := \sum_{k=0}^K c_k \cos(2\pi k \tau) \mathbf 1_{[-1/2,1/2)}.
  \end{equation}
  Then for any integer $R > 2K$
  \begin{equation}\label{eq:sumSquares}
    \sum_{m\in {\mathbb Z}} \left|\vartheta\left(\frac{\tau-m}{R}\right)\right|^2 = R c_0^2 + \frac{R}{2} \sum_{k=1}^K c_k^2, ~~\forall \tau \in \RR;
  \end{equation}
  i.e., the sum of squares of a system of regular translates $\left( \bd T_m \theta \right)_{m\in \ZZ}$, with $\theta = \vartheta(\cdot/R)$, is constant.

\revised{Considering Theorem \ref{thm:PLwarp}(iii), the construction above can be used to} easily construct tight frames by choosing the decimation factors $a_m$ to satisfy
\begin{equation}\label{eq:maxam}
  a_m^{-1} \geq \Phi^{-1}(m+R/2)-\Phi^{-1}(m-R/2).
\end{equation}
In the following, we will demonstrate this for some of the examples given in Section \ref{sec:warp}. 

For the purpose of all the following examples, we choose $\vartheta$ according to \eqref{eq:window_form} with $K=1$ and $c_0=c_1=1/2$, i.e., we can choose $R\geq 3$.
This function is often called the \emph{Hann} or raised cosine window. The Hann window is among the most popular finitely supported Gabor windows or filters for time-frequency signal analysis.

\begin{example}[$\Phi(\xi) = \sgn(\xi)\log(1+|\xi|)$]\label{ex:ERBtight} For this choice of $\Phi$, \eqref{eq:maxam} takes the form
\begin{equation*}
    a_m^{-1} \geq \sgn(m+R/2)(e^{|m+R/2|}-1)-\sgn(m-R/2)(e^{|m-R/2|}-1), 
\end{equation*}
or equivalently
\begin{equation*}
    a_m^{-1} \geq \begin{cases}
      (e^{|m|+R/2}-1) - (e^{|m|-R/2}-1) = e^{|m|}(e^{R/2}-e^{-R/2}) & \text{ for } |m|\geq R/2,\\
      (e^{m+R/2}-1) + (e^{-m+R/2}-1) = e^{R/2}(e^{|m|}+e^{-|m|})-2 & \text{ else,}
    \end{cases}
\end{equation*}
where the latter case concerns the filters where $\supp(\bd T_m\theta)$ is not contained in either $[0,\infty)$ or $(-\infty,0]$. We see that in both cases, $a_m^{-1}$ is majorized by $e^{|m|}$, up to a constant depending solely on $R$. If we set $R=3$, then a tight frame is obtained by choosing 
  \begin{equation*}
    a_m = \begin{cases}
            e^{-|m|}(e^{3/2}-e^{-3/2})^{-1} \geq \frac{1}{4.26 e^{|m|}} & \text{ for } |m|\geq 2,\\
            (e^{3/2}(e^{1}+e^{-1})-2)^{-1} > \frac{1}{11.84} & \text{ for } |m|=1,\\
            (2e^{3/2}-2)^{-1} > \frac{1}{6.97} & \text{ for } m=0.
          \end{cases}
  \end{equation*}
\end{example}
On the other hand, Corollary \ref{cor:PLnatural} yields $\widetilde{a_w}=(2e^{3/2}-2)^{-1}$ and $w(m)=v(m) = e^{|m|}$, i.e., $a_m = 2e^{-|m|}(e^{3/2}-1)^{-1}$ for all $m\in\ZZ$, which is slightly more conservative.

\begin{example}[$\Phi_\alpha(\xi) = \sgn(\xi)((1+|\xi|)^{1-\alpha}-1)$]
  Let $p:=1/(1-\alpha)\in\NN$. Then  \eqref{eq:maxam} can be rewritten as 
  \begin{equation*}    
    a_m^{-1} \geq \begin{cases} 
     (1+|m|+R/2)^p - (1+|m|-R/2|)^p & \text{ for } |m|\geq R/2,\\ 
     (1+R/2+m)^p + (1+R/2-m)^p -2 & \text{ else.}
     \end{cases}
  \end{equation*}
  If $\alpha=1/2$, i.e., $p=2$, and $R=3$, evaluation of the above conditions yields a tight warped filter bank with 
   \begin{equation*}
    a_m = \begin{cases}
            \frac{1}{6+6|m|} & \text{ for } |m|\geq 2,\\
            \frac{2}{25} & \text{ for } |m|=1,\\
            \frac{2}{21} & \text{ for } m=0.
          \end{cases}
  \end{equation*}
  In this setting, Corollary \ref{cor:PLnatural} yields $\widetilde{a_w}=\frac{4}{21}$ and $w(m) = 2+2|m|$, see also Example \ref{ex:alpha1}. 
\end{example}


\begin{example}[$\Phi(\xi) = \tan(\xi)$] 
  With $w$ and $v$ as in Proposition \ref{pro:warpingfunctions}(v), Corollary \ref{cor:PLnatural} yields 
  $\widetilde{a_w} = (R^3/6 + 2R - 4)^{-1}$ and, with $R=3$, the following set of almost optimal natural decimation factors
  \[
   a_m = \widetilde{a_w}/w(m) = \frac{2+2m^2}{13}.
  \]
 \end{example}
 
 The above examples show the ease with which tight warped filter bank frames are constructed when the prototype filter $\theta$ is compactly supported. In the following section, we show that warped filter bank frames for fully supported prototype filters exist as well.  
 
 \subsection{Prototype decay implies existence of warped filter bank frames}\label{ssec:framesFromDecay}   

Our final set of sufficient Bessel and frame conditions is concerned with the case that $\theta$ sufficiently localized, but not necessarily compactly supported. In this setting, the verification of the frame property becomes substantially harder. To obtain a sufficient condition, it is possible to estimate the alias terms $\sum_{l \neq 0} \abs{g_m\overline{\bd T_{la_m^{-1}} g_m}}$, $m\in\ZZ$, in the Walnut representation of the frame operator of $\mathcal G(\Phi,\theta,\bd a)$, see \cite[Proposition 3.7]{jale14} \revised{a variant of which we state in Proposition \ref{pro:jaleProp}}. The main result of this section, Theorem \ref{thm:wpdsuffcond} provides a decay condition on $\theta$ and a density condition on the decimation factors $\bd a$, such that the conditions of \cite[Proposition 3.7]{jale14} are satisfied. Note that these conditions were recently improved by Lemvig et al.~\cite{leve17}, under the additional assumption of the so-called $\alpha$-local integrability 
condition. However, our results do not benefit from the sharper condition.

One would be tempted to apply~\cite[Corollary 3.5]{doma14-1}, which uses decay of the frequency responses $g_m$ to determine a density condition on the decimation factors $a_m$. Adapted to our setting, that result is as follows:

\begin{proposition}[\cite{doma14-1}, Corollary 3.5]\label{pro:polybound}
  Let $(g_m)_{m\in\ZZ}\subset \LD$ satisfying 
  \begin{equation}\label{eq:gmboundnoam}
   0 < \tilde{A} \leq \sum_{m\in\ZZ} |g_m(\xi)|^2 \leq \tilde{B} <\infty \quad \text{a.e. on } D,
  \end{equation}
  for some constants $\tilde{A},\tilde{B}$. Assume that there exist a $\delta$-separated sequence $(b_m)_{m\in\ZZ}$ and constants $p>2$ and $C>0$, such that 
  \begin{equation}\label{eq:polybound}
   |g_m(\xi)| \leq C (1 + |\xi-b_m|)^{-p}\quad \text{a.e. on } D.
  \end{equation}
  There exists a sequence $(a^{(0)}_m)_{m\in\ZZ}\subset \RR^+$ such that the filter bank $(g_{m,n})_{m,n\in\ZZ}$ generated from $(g_m)_{m\in\ZZ}\subset \LD$ and $(a_m)_{m\in\ZZ}\subset \RR^+$ forms a frame for $\LtDF$, if $a_m \leq a^{(0)}_m$, for all $m\in\ZZ$.
\end{proposition}

Note that we took the liberty of simplifying the conditions of \cite[Corollary 3.5]{doma14-1}: In particular, the assumption that the $g_m$ be in the Wiener space $W(\bd L^{\infty},\ell^1)$ is implied by \eqref{eq:polybound}. Moreover, the constants $C$ and $p$ were allowed to vary with $m\in\ZZ$, but only within a compact interval. It is straightforward to confirm that the two results are equivalent if $D = \RR$, except for the actual values of the sequence $(a^{(0)}_m)_{m\in\ZZ}\subset \RR^+$. The required restriction in the case $D\subsetneq \RR$ is straightforward as well. 

The conditions of Proposition \ref{pro:polybound} pose severe restrictions for warped filter banks. In fact, under reasonable assumptions on $\mathcal G(\Phi,\theta,\bd a)$, they imply that $w=(\Phi^{-1})'$ is bounded above and thus $\Phi$ must have at least linear asymptotic growth. 

\begin{proposition}\label{pro:linGrowth}
  Let $\mathcal G(\Phi,\theta,\bd a)$ be a warped filter bank for $\LtDF$, with nonzero prototype $\theta$, and set $b_m = \Phi^{-1}(m)$, for all $m\in\ZZ$. 
  \begin{itemize}
   \item[(i)] If the open interval $D$ is a true subset of $\RR$, i.e., $D\subsetneq \RR$, then for any $\delta>0$, there is an $m\in\ZZ$, such that $|b_{m+1}-b_m|\leq \delta$.
   \item[(ii)] Assume that $\theta \leq C_0 (1+|\cdot|)^{-p_0}$, for $C_0>0$, $p_0>1/2$, and that  \eqref{eq:gmboundnoam} holds for $g_m = \sqrt{a_m}(\bd T_m \theta)\circ \Phi$, $m\in\ZZ$. Then $\sup_{m\in\ZZ} a_m < \infty$ and $\limsup_{m\rightarrow \infty} a_m \neq 0 \neq \limsup_{m\rightarrow -\infty} a_m$.
   \item[(iii)] If the assumptions of (ii) hold and there are constants $p,C>0$, such that $|g_m| \leq C(1+|(\cdot)-b_m|)^{-p}$ almost everywhere, for all $m\in\ZZ$, then $w = (\Phi^{-1})' \in \bd L^\infty(\RR)$.
  \end{itemize}
\end{proposition}
\begin{proof}
  Ad (i): Assume without loss of generality that $D$ is bounded below with $\inf\{\xi: \xi\in D\} = c \in \RR$. Clearly, since $\Phi'$ is is continuous and positive, $\lim_{\xi\rightarrow c} \Phi'(\xi) = \infty$, implying $\lim_{\tau\rightarrow -\infty} (\Phi^{-1})'(\tau) = 0$ and (i) easily follows.
  
  Ad (ii): Without loss of generality, assume $\esssup_{\tau\in\RR} \theta(\tau) = 1$. Then $\sup_{m\in\ZZ} a_m = \infty$ or $\lim_{m\rightarrow -\infty} a_m = \infty$ easily imply that a finite upper bound for $\sum_{m\in\ZZ} |g_m|^2 = \sum_{m\in\ZZ} a_m|\theta(\Phi(\cdot) - m)|^2$ cannot exist. Let $\overline{a} := \sup_{m\in\ZZ} a_m < \infty$. For every $\epsilon > 0$, there is $m_{\epsilon}\in\ZZ$, such that 
  \[
   \sum_{m\leq (k-m_{\epsilon})} a_m|\theta(\tau - m)|^2 \leq C_0^2\sum_{m\leq (k-m_{\epsilon})} \overline{a}(1+|\tau - m|)^{-2p_0} \leq \epsilon/2, \text{ for almost every } \tau\geq k,\ k\in\ZZ.
  \]
 If $\limsup_{m\rightarrow \infty} a_m = 0$, then there is a $k_\epsilon\in\ZZ$, such that $\sup\limits_{m>(k_\epsilon-m_{\epsilon})} a_m  \leq \epsilon/(2B)$, where the constant $B$ is defined as
  \[
   B := \esssup_{\tau\in\RR} \sum_{m\in\ZZ} |\theta(\tau-m)|^2 \leq \esssup_{\tau\in\RR} C_0^2\sum_{m\in\ZZ} (1+|\tau-m|)^{-2p_0} < \infty.
  \]
  Together, we obtain 
  \[
   \sum_{m\in\ZZ} a_m|\theta(\tau - m)|^2 \leq \epsilon, \text{ for almost every } \tau \geq k_\epsilon.
  \]
  Since $\epsilon>0$ is arbitrary, the desired lower bound $\tilde{A}$ cannot exist. $\limsup_{m\rightarrow -\infty} a_m \neq 0$ is proven by the same steps.
  
  Ad (iii): We show that under the assumptions of (iii), $w\notin \bd L^\infty(\RR)$ implies $\theta \equiv 0$, which contradicts the assumption that $\theta$ is nonzero. Begin by noting that $|g_m| \leq C(1+|(\cdot)-\Phi^{-1}(m)|)^{-p}$ is equivalent to 
  \[
   \theta \leq \frac{C}{\sqrt{a_m}(1+|\Phi^{-1}(\cdot+m)-\Phi^{-1}(m)|)^{p}}, \text{ almost everywhere}.
  \]
  Moreover, for all $\tau_0\in\RR$, 
  \[
   \Phi^{-1}(\tau_0+m)-\Phi^{-1}(m) = \int_{0}^{\tau_0} w(m+\tau)~ d\tau \geq \int_0^{\tau_0} \frac{w(m)}{v(-\tau)}~ d\tau,
  \]
  where we used $v$-moderateness of $w$. If $\tau \geq \tau_1 > 0$, then with $C_{\tau_1} := \int_0^{\tau_1} \frac{1}{v(-\tau)} d\tau$, we obtain by positivity of $v$ that 
  \begin{equation}\label{eq:boundfortheta}
   \theta(\tau) \leq \frac{C}{\sqrt{a_m}(1+|\Phi^{-1}(\tau+m)-\Phi^{-1}(m)|)^{p}} \leq \frac{C}{\sqrt{a_m}(1+C_{\tau_1}w(m)|)^{p}}, \text{ for almost all } \tau\geq \tau_1.
  \end{equation}
  Now, if $w\notin \bd L^\infty(\RR)$, then either $\lim_{m\rightarrow \infty} w(m) = \infty$ or $\lim_{m\rightarrow -\infty} w(m) = \infty$, by continuity of $w$. For the right hand side of  \eqref{eq:boundfortheta}  to be bounded below, this implies either $\lim_{m\rightarrow \infty} a_m = 0$ or $\lim_{m\rightarrow -\infty} a_m = 0$, which is prohibited by (ii). Since $\tau_1>0$ was arbitrary, we obtain that necessarily $\theta(\tau) = 0$ for almost every $\tau>0$. An analogous argument shows that $\theta(\tau) = 0$ for almost every $\tau<0$. Therefore, $\theta\equiv 0$, completing the proof by contradiction.
\end{proof}

Considering Proposition \ref{pro:linGrowth}, a quick glance at Examples \ref{ex:wavelet1}--\ref{ex:alpha1} shows that the requirements of Proposition \ref{pro:polybound} are highly undesirable for warped filter banks. \revised{Instead, we only establish a decay condition on $\theta$.} This condition ensures the Bessel property and, when complemented by sufficiently small decimation factors, even the frame property. The given result is of central interest, as it shows that warped filter banks also admit the construction of frames for many prototype filters $\theta$ with full support. 

\begin{theorem}\label{thm:wpdsuffcond}
  Let $\Phi: D\rightarrow \RR$ be a warping function with $\theta\in \bd L^2_{\sqrt{v}}(\RR)$, fix an arbitrary $\epsilon>0$ and let $w_1,w_2$ denote the following weights
  \[
   w_1 = (1+|\cdot|)^{1+\epsilon} \quad \text{ and }\quad w_2 = (1+|V(0,\cdot)|)^{1+\epsilon},
  \]
  where $V$ is as defined in  \eqref{eq:defOfV}. If 
  \begin{equation}\label{eq:DecayAndNatFacs} 
    \theta\in \bd L^\infty_{w_1}(\RR)\cap \bd L^\infty_{w_2}(\RR)\quad \text{and} \quad a_m \leq \tilde a/w(m),\quad \text{ for all } m\in\ZZ \text{ and some }\tilde a>0,
  \end{equation}
  then $\mathcal{G}(\Phi,\theta,\bd{a})$ is a Bessel sequence. If \eqref{eq:DecayAndNatFacs} holds and  
  \[
    0 < A_1 \leq \sum_{m\in\ZZ}  |\bd T_m \theta|^2 \text{ almost everywhere,}
  \]
  then there is a constant $\tilde a_{0}>0$ such that $\mathcal{G}(\Phi,\theta,\bd{a})$ is a frame, whenever $a_m\leq \tilde a_{0}/w(m)$, for all $m\in\ZZ$.
\end{theorem}

\revised{We would like to emphasize that \eqref{eq:DecayAndNatFacs} is only a mild additional restriction on $\theta$, given that $\theta\in\mathbf L^2_{\sqrt{v}}(\RR)$ is already required. Clearly, $\theta \in \bd L^\infty_{w_1}(\RR)$ enforces just enough decay to ensure $\theta\in \bd L^1(\RR)$, while simple calculations show that, similarly, $\theta \in \bd L^\infty_{w_2}(\RR)$ provides just enough decay to ensure $\theta\in \bd L_v^1(\RR)$.}


Before we proceed to prove Theorem \ref{thm:wpdsuffcond}, we require two auxiliary results.

\begin{lemma}\label{lem:PhiDistEst}
  Let $\Phi:D\rightarrow \RR$ be a warping function such that $w$ is $v$-moderate. There is a bijective, increasing function $A_v:\RR\rightarrow \RR$, such that $A_v(0) = 0$ and for all $c\in\RR^+$, we have
   \[\left|\Phi^{-1}(\tau_1)-\Phi^{-1}(\tau_0)\right|\geq cw(\tau_0)\quad \Longrightarrow \left|\tau_1-\tau_0\right| \geq \begin{cases}
      |A_v^{-1}(c)| & \text{ if } \tau_1 \geq \tau_0\\
      |A_v^{-1}(-c)| & \text{ else}
   \end{cases}, \text{ for all } \tau_0,\tau_1\in\RR.\]
\end{lemma}
\begin{proof}
  If $\tau_1\geq\tau_0$, then the assumptions yield
  \[
   cw(\tau_0) \leq \Phi^{-1}(\tau_1) - \Phi^{-1}(\tau_0) = \int_{\tau_0}^{\tau_1} w(\tau) d\tau \leq w(\tau_0)(\tau_1 - \tau_0) \sup_{\tau\in[0,\tau_1-\tau_0]} v(\tau).
  \]
  Analogous, we obtain for $\tau_1 < \tau_0$ that $cw(\tau_0) \leq w(\tau_0)(\tau_0 - \tau_1) \sup_{\tau\in[\tau_1-\tau_0,0]} v(\tau)$. 
  
  The function $A_v : \RR \rightarrow \RR$, $\tau\mapsto \tau \sup_{\tau_0 \in[\tau,0]\cup [0,\tau]} v(\tau_0)$ is continuous and strictly increasing and thus invertible. Moreover, the above derivations show that, for all $\tau_0,\tau_1\in\RR$,
  \[
   c \leq \sgn(\tau_1-\tau_0)A_v(\tau_1-\tau_0),
  \]
  as desired. 
\end{proof}

Lemma \ref{lem:PhiDistEst} allows us to derive the following result which will be crucial for proving Theorem \ref{thm:wpdsuffcond}. 

\begin{lemma}\label{lem:offdiagest}
   For a given warped filter bank $\mathcal G(\Phi,\theta,\bd a)$, with a sequence $\bd a = (a_m)_{m\in\ZZ}$ of decimation factors, define
   \begin{equation}\label{eq:defOfP}
      \bd P(\xi) := \bd P_{\Phi,\theta,\bd a}(\xi) := \sum_{m\in\ZZ} \left(\left|\theta(\Phi(\xi)-m)\right|\cdot \sum_{\substack{ k\in\ZZ\setminus\{0\}\\ \xi+ka_m^{-1}\in D}} \left|\theta(\Phi(\xi+ka_m^{-1})-m)\right|\right), \text{ for all } \xi\in D.
    \end{equation}
    If $\theta\in \bd L^\infty_{w_1}(\RR)\cap \bd L^\infty_{w_2}(\RR)$, with $w_1,w_2$ as in Theorem \ref{thm:wpdsuffcond}, and $a_m \leq \tilde a/w(m)$, for all $m\in\ZZ$ and some $\tilde a>0$, then 
    \[
      \esssup_{\xi\in D} \bd P(\xi) < \infty\quad \text{ and }\quad  \esssup_{\xi\in D} \bd P(\xi) \overset{\tilde{a}\rightarrow 0}{\longrightarrow} 0.
    \]
\end{lemma}
\begin{proof}
  In the derivations below, sums of the kind $\sum_{k\in\NN_0} (q+k)^{-s}$, for $q>0,s>1$ appear repeatedly. These sums are finite and their value is given by \emph{Hurwitz' zeta function}~\cite{ol10-2}, $\zeta(q,s)$.
  In the following, we will use the upper estimate
  \begin{equation}
    q^{-s} + \sum_{k\in\NN} (q+k)^{-s} < q^{-s} + \int_{\RR^+} (q+t)^{-s}~dt = q^{-s} + (s-1)^{-1}q^{1-s},
    \label{eq:HurwitzEst}
  \end{equation}
  which also shows that the left hand side tends towards zero for $s$ fixed and $q\rightarrow \infty$ or vice versa. Note that $v$-moderateness of $w$ implies, for $t\in D$,  
  \begin{align}
     \left|t-\Phi^{-1}(m)\right| & = \left|\Phi^{-1}(\Phi(t))-\Phi^{-1}(m)\right| 
     = \left|\int_m^{\Phi(t)} w(\tau)~d\tau \right|\nonumber\\
     & = \left|\int_0^{\Phi(t)-m} w(\tau+m)~d\tau \right|\nonumber\\
     & \leq w(m)\left|\int_0^{\Phi(t)-m} v(\tau)~d\tau\right| = w(m)|V(0,\Phi(t)-m)|.
    \label{eq:innerEst1}
  \end{align}
  Moreover, with $C_1 := \|\theta\|_{\bd L^\infty_{w_1}(\RR)} < \infty$,
  \begin{equation}\label{eq:decupperbound}
    \esssup_{\tau\in\RR} \sum_{m\in\ZZ} |\bd T_m \theta(\tau)| \leq 2C_1\sum_{m\in\NN_0} \frac{1}{(1+m)^{1+\epsilon}} < 2C_1(1+\epsilon^{-1}) =: \tilde{B} < \infty,
  \end{equation}  
  by \eqref{eq:HurwitzEst}.
  We proceed to estimate the inner sum in \eqref{eq:defOfP}. 
  To that end, define
  \[
   P_m(\xi):= \sum_{\substack{ k\in\ZZ\setminus\{0\}\\ \xi+ka_m^{-1}\in D}} \left|\theta(\Phi(\xi+ka_m^{-1})-m)\right|, \text{ for all } \xi\in D,\ m\in\ZZ.
  \]
  Assuming $\theta\in \bd L^\infty_{w_2}(\RR)$ with $C_2: = \|\theta\|_{\bd L^\infty_{w_2}(\RR)}  > 0$, we obtain
  \[
   \begin{split}
   P_m(\xi) & \leq C_2 \sum_{\substack{ k\in\ZZ\setminus\{0\}\\ \xi+ka_m^{-1}\in D}} \left|(1+ |V(0,\Phi(\xi+ka_m^{-1})-m))|)^{-1-\epsilon}\right|\\
   & \leq C_2 \sum_{k\in\ZZ\setminus\{0\}} \left(1+ \left|\frac{\xi+ka_m^{-1}-\Phi^{-1}(m)}{w(m)}\right|\right)^{-1-\epsilon},
   \end{split}
  \]
  for almost every $\xi\in D$. Here, we used \eqref{eq:innerEst1} with $t=\xi+ka_m^{-1}$ to obtain the second inequality. 
  
  For any pair $(\xi,m)$, there is a unique $k_{(\xi,m)}\in\ZZ$ such that $\xi + k_{(\xi,m)}a_m^{-1}-\Phi^{-1}(m) \in [-(2a_m)^{-1},(2a_m)^{-1})$. Let 
  \[
   M_\xi := \{m\in\ZZ \colon  k_{(\xi,m)} = 0\}\quad \text{ and }\quad M_\xi^\dag := \ZZ\setminus M_\xi.
  \]
  First assume that $m\in M_\xi$, i.e., $\xi-\Phi^{-1}(m)\in [-(2a_m)^{-1},(2a_m)^{-1})$. We can split the sum by the sign of $k$ to obtain
  \[
  \begin{split}
   \lefteqn{\sum_{k\in\ZZ\setminus\{0\}} \left(1+ \left|\frac{\xi+ka_m^{-1}-\Phi^{-1}(m)}{w(m)}\right|\right)^{-1-\epsilon}}\\
   & \leq \sum_{k\in\NN_0} \left(1+\left|\frac{(2a_m)^{-1}+ka_m^{-1}}{w(m)}\right|\right)^{-1-\epsilon} + \sum_{k\in\NN_0} \left(1+ \left|-\frac{(2a_m)^{-1}+ka_m^{-1}}{w(m)}\right|\right)^{-1-\epsilon}\\
   & = 2 \sum_{k\in\NN_0} \left(\frac{w(m)+(2a_m)^{-1}+ka_m^{-1}}{w(m)}\right)^{-1-\epsilon} = (\ast).
   \end{split}
  \]
  If $a_m \leq \tilde a /w(m)$, then $|(2a_m)^{-1}+ka_m^{-1}| \geq |(1/2+k)w(m)/\tilde a|$ and 
  \[
   (\ast) \leq 2 \sum_{k\in\NN_0} \left(\frac{\tilde{a}+1/2+k}{\tilde{a}}\right)^{-1-\epsilon} = 2 \tilde{a}^{1+\epsilon}\sum_{k\in\NN_0} \left(\tilde{a}+1/2+k\right)^{-1-\epsilon}.
  \]
  Consequently, with \eqref{eq:HurwitzEst}, we obtain that 
  \[
   P_m(\xi) 
   < 2C_2 \left(1+\frac{\tilde{a}+1/2}{\epsilon}\right)\left(\frac{\tilde{a}}{\tilde{a}+1/2}\right)^{1+\epsilon}, 
  \]
  for almost every $\xi\in D$. 
  Now, if $m\in M_\xi^\dag$, then a similar estimation yields
   \[
   \begin{split}
   P_m(\xi) & \leq \left|\theta(\Phi(\xi+k_{(\xi,m)}a_m^{-1})-m)\right| + C_2\sum_{k\in\ZZ\setminus\{0,k_{(\xi,m)}\}} \left(1+ \left|\frac{\xi+ka_m^{-1}-\Phi^{-1}(m)}{w(m)}\right|\right)^{-1-\epsilon}\\
   & 
   < C_2 + 2C_2 \left(1+\frac{\tilde{a}+1/2}{\epsilon}\right)\left(\frac{\tilde{a}}{\tilde{a}+1/2}\right)^{1+\epsilon}, 
   \end{split}
  \]
  almost everywhere.  
  These estimates can now be inserted into the expression \eqref{eq:defOfP} for $\bd P$:
  \begin{equation}\label{eq:PEstimate}
   \begin{split}
     \bd P(\xi) & = \sum_{m\in\ZZ} |\theta(\Phi(\xi)-m)|\cdot P_m(\xi)\\
       & \leq C_2 \sum_{m\in M_\xi^\dag} |\theta(\Phi(\xi)-m)| + 2C_2 \left(1+\frac{\tilde{a}+1/2}{\epsilon}\right)\left(\frac{\tilde{a}}{\tilde{a}+1/2}\right)^{1+\epsilon}\cdot \sum_{m\in\ZZ} |\theta(\Phi(\xi)-m)|,
   \end{split}  
  \end{equation}
  for almost every $\xi\in D$. 
  Applying \eqref{eq:decupperbound} to estimate both terms yields 
  \[
   \begin{split}
   \esssup_{\xi\in D} \bd P(\xi) 
   & < 2 C_1 C_2 (1+\epsilon^{-1})\left(1+ 2\left(1+\frac{\tilde{a}+1/2}{\epsilon}\right)\left(\frac{\tilde{a}}{\tilde{a}+1/2}\right)^{1+\epsilon}\right) < \infty.
   \end{split}
  \]
  For the second assertion, we have to show convergence of the essential supremum to $0$ for $\tilde{a}\rightarrow 0$. By definition, $m\in M_\xi^\dag$ implies $|\xi - \Phi^{-1}(m)| \geq (2a_m)^{-1} \geq w(m)(2\tilde a)^{-1}$. Hence, we can apply Lemma \ref{lem:PhiDistEst}, with $\tau_0 = m,\ \tau_1 = \Phi(\xi)$ and $c=(2\tilde a)^{-1}$ to obtain
  \[
   \left|\Phi(\xi)-m\right| \geq \begin{cases}
                                   |A_v^{-1}\left(1/(2\tilde{a})\right)| & \text{ if } \Phi(\xi)-m \geq 0,\\
                                   |A_v^{-1}\left(-1/(2\tilde{a})\right)| & \text{ else. }
                                 \end{cases}
  \]
  We can rewrite, with $m_- := \max M_\xi^\dag\cap (-\infty,\Phi(\xi)), m_+ := \min M_\xi^\dag\cap (\Phi(\xi),\infty)$, 
  \[
   \begin{split}
   \sum_{m\in M_\xi^\dag} |\theta(\Phi(\xi)-m)| & \leq \sum_{k\in\NN_0} |\theta(\Phi(\xi)-m_- +k)| + \sum_{k\in\NN_0} |\theta(\Phi(\xi)-m_+ -k)|\\   
   & \leq C_1\left(\sum_{k\in\NN_0} \frac{1}{\left(1+|A_v^{-1}\left(1/(2\tilde{a})\right)|+k\right)^{1+\epsilon}} + \sum_{k\in\NN_0} \frac{1}{\left(1+|A_v^{-1}\left(-1/(2\tilde{a})\right)|+k\right)^{1+\epsilon}}\right)\\
   & < C_1 \cdot \sum_{j=0}^1 \left(1 + \frac{1+|A_v^{-1}\left((-1)^{j}/(2\tilde{a})\right)|}{\epsilon}\right)\left(1+|A_v^{-1}\left((-1)^{j}/(2\tilde{a})\right)|\right)^{-1-\epsilon}
   \end{split}
  \]
  All in all, we obtain for $\bd P$ the estimate
  \begin{equation}
    \begin{split}
    \lefteqn{\bd P_{\Phi,\theta,\bd a}(\xi)}\\
   & < C_1C_2 \left(4(1+\epsilon^{-1})\left(1+\frac{\tilde{a}+1/2}{\epsilon}\right)\left(\frac{\tilde{a}}{\tilde{a}+1/2}\right)^{1+\epsilon} + \sum_{j=0}^1 \frac{1+\epsilon^{-1}(1+|A_v^{-1}\left((-1)^{j}/(2\tilde{a})\right)|)}{(1+|A_v^{-1}\left((-1)^{j}/(2\tilde{a})\right)|)^{1+\epsilon}}\right),
    \end{split}\label{eq:FinalPEstimate}
  \end{equation}
  almost everywhere.
  Since $A_v^{-1}(\tau) \overset{\tau\rightarrow \pm\infty}{\longrightarrow} \infty$, we see that 
  \[
   \esssup_{\xi\in D} \bd P_{\Phi,\theta,\bd a}(\xi) \overset{\tilde{a}\rightarrow 0}{\longrightarrow} 0,
  \]
  as desired, finishing the proof.
\end{proof}

The last term on the right hand side of  \eqref{eq:FinalPEstimate} depends heavily on the moderating weight $v$ through the function $A_v^{-1}$ and without further specifying $v$, a useful estimate for $A_v^{-1}$ is out of reach. For $v(\tau) = e^{\tau}$, cf. Example \ref{ex:wavelet1}, we have $A_v(\tau) = \tau \max\{1,e^\tau\}$. Thus, on $\RR^+$, $A_v^{-1}$ equals the product logarithm, such that the right hand side of  \eqref{eq:FinalPEstimate} will decay very slowly for $\tilde{a}\rightarrow 0$. \revised{However, it decays quickly for increasing $\epsilon$. In our experiments, performed with compactly supported or exponentially decaying $\theta$, we never observed any significant influence of the warping function $\Phi$ (and therefore of $v$) on the choice of $\tilde{a}$. These observations are reflected in the results presented in Section \ref{ssec:ExpRatios}, in particular the frame bound ratios reported in Table \ref{tab:framebounds}.} On the other hand, the estimate \eqref{eq:FinalPEstimate} may be 
rather coarse. For a smooth bell function $\theta$, e.g. a Gaussian, even the base estimates $\theta \leq C_0(1+|\cdot|)^{-1-\epsilon}$ and $\theta \leq C_1(1+|V(0,\cdot)|)^{-1-\epsilon}$ do not allow the simultaneous choice of small constants $C_0,C_1$ and a large decay rate $\epsilon$.

With Lemma \ref{lem:offdiagest} in place, proving Theorem \ref{thm:wpdsuffcond} only requires a few simple steps and the application of the following variant of \cite[Proposition 3.7]{jale14}.

\begin{proposition}[\cite{jale14}, Proposition 3.7]\label{pro:jaleProp}
  Let $(g_{m,n})_{m,n\in\ZZ}\subset \LtDF$ be the filter bank generated from $(g_m)_{m\in\ZZ}\subset \LD$ and $(a_m)_{m\in\ZZ}\subset \RR^+$. If  
  \begin{equation}\label{eq:suffupper}     
      B := \esssup_{\xi\in D}\left[ \sum_{m\in\ZZ}\sum_{l \in\ZZ} \frac{1}{a_m}\abs{g_m(\xi)\overline{g_m}(\xi-l/a_m)}\right] < \infty,
 \end{equation}
 then $(g_{m,n})_{m,n\in\ZZ}$ is a Bessel sequence with bound $B$. 
 
 Assume that \eqref{eq:suffupper} holds. If 
 \begin{equation}\label{eq:sufflower}     
    A := \essinf_{\xi\in D}\left[ \sum_{m\in\ZZ}\frac{1}{a_m}\left(\abs{g_m(\xi)}^2-\sum_{l\in\ZZ\setminus\{0\}} \abs{g_m(\xi)\overline{g_m}(\xi-l/a_m)} \right)\right] > 0,
 \end{equation}
 then $(g_{m,n})_{m,n\in\ZZ}$ constitutes a frame for $\LtDF$ with frame bounds $A,B$.
\end{proposition}

\begin{proof}[Proof of Theorem \ref{thm:wpdsuffcond}]
 We will show that, with suitable choices of $\tilde{a}$, the conditions of 
 Theorem \ref{thm:wpdsuffcond} enable the application of Proposition \ref{pro:jaleProp}.  
  The main observation is the following: For any given warped filter bank $\mathcal G(\Phi,\theta,\bd a)$, we have 
  \[
   \begin{split}
   \lefteqn{\sum_{m\in\ZZ} a_m^{-1}|g_m(\xi)|^2 \pm \sum_{m\in\ZZ}\sum_{l\in\ZZ\setminus\{0\}} a_m^{-1}|g_m(\xi)\overline{g_m(\xi+la_m^{-1})}|}\\
  & = \sum_{m\in\ZZ} |\theta(\Phi(\xi)-m)|^2 \pm \sum_{m\in\ZZ} \left(|\theta(\Phi(\xi)-m)| \cdot \sum_{\substack{ l\in\ZZ\setminus\{0\}\\ \xi+la_m^{-1}\in D}} |\theta(\Phi(\xi+la_m^{-1})-m)|\right)\\
  & = \sum_{m\in\ZZ} |\theta(\Phi(\xi)-m)|^2 \pm \bd P_{\Phi,\theta,\bd a}(\xi), \text{ for almost every } \xi\in D.
   \end{split}
  \]
  By Lemma \ref{lem:offdiagest}, $\esssup_{\xi\in D} \bd P_{\Phi,\theta,\bd a}(\xi) < \infty$. Moreover, since $\theta\in\bd L^\infty_{w_1}(\RR)$, we obtain the estimate
  \[
    \esssup_{\tau\in\RR} \sum_{m\in\ZZ} |\bd T_m \theta(\tau)| \leq \tilde{B} < \infty
  \]
  as per \eqref{eq:decupperbound}. In total, with $B$ as in \eqref{eq:suffupper},
  \[
    B = \esssup_{\xi\in D}\left( \sum_{m\in\ZZ} |\theta(\Phi(\xi)-m)|^2 + \bd P_{\Phi,\theta,\bd a}(\xi)\right) \leq \tilde{B}^2 + \esssup_{\xi\in D} \bd P_{\Phi,\theta,\bd a}(\xi) < \infty,
  \]
  and $\mathcal G(\Phi,\theta,\bd a)$ is a Bessel sequence by Proposition \ref{pro:jaleProp}. Similarly, with $A$ as in \eqref{eq:sufflower}, 
  \[
   A = \essinf_{\xi\in D}\left( \sum_{m\in\ZZ} |\theta(\Phi(\xi)-m)|^2 - \bd P_{\Phi,\theta,\bd a}(\xi)\right)\\ \geq A_1 - \esssup_{\xi\in D} \bd P_{\Phi,\theta,\bd a}(\xi).
  \]
  By Lemma \ref{lem:offdiagest}, there is a constant $\tilde{a}_{0}>0$, such that $a_m\leq \tilde{a}_{0}/w(m)$, for all $m\in\ZZ$, implies
  \[
    \esssup_{\xi\in D} \bd P_{\Phi,\theta,\bd a}(\xi) < A_1.
  \]
  Thus, by Proposition \ref{pro:jaleProp}, we have that $\mathcal G(\Phi,\theta,\bd a)$ constitutes a frame. 
\end{proof}

Theorem \ref{thm:wpdsuffcond} is extremely useful for proving (a) the \emph{existence of a safe region} in which decimation factors can be chosen freely and (b) that compact support of the prototype $\theta$ is not a necessity for obtaining warped filter bank frames.

\subsection{Warped filter bank frames with $\theta\in \bd L^2_{\sqrt{w}}(\RR)$}\label{sec:ThetaInLtw}

When going through the results presented in this section, we regularly use the moderateness of $w=(\Phi^{-1})'$ to obtain estimates $w(\tau+m) \leq w(m)v(\tau)$. Clearly, we can exchange the roles of $w$ and $v$, to obtain estimates in terms of $v(m)$ instead of $w(m)$. With this simple change, we can recover all the presented results in the setting where $\theta\in \bd L^2_{\sqrt{w}}(\RR)$ and natural decimation factors take the form $a_m = \tilde{a}/v(m)$, for some $\tilde{a}>0$. Adapting the proofs amounts to simply exchanging the roles of $w$ and $v$. Furthermore, $V(\tau_0,\tau_1)$ must be exchanged for $W(\tau_0,\tau_1):= \Phi^{-1}(\tau_1)-\Phi^{-1}(\tau_0)$, in the statements and proofs of Corollary \ref{cor:PLnatural}, Lemma \ref{lem:offdiagest} and Theorem \ref{thm:wpdsuffcond}.

\section{Warped filter banks for digital signals}\label{sec:discWarp}

We now turn our attention to discrete signals, i.e., sequences $x\in\ell^2(\ZZ)$, the customary setting for studying filter banks. The discrete time Fourier transform (DTFT) $\hat x(\xi) := \mathcal{F}_d x(\xi)
= \sum_\ZZ x(l) e^{-2 \pi i l \xi/\xi_s}$ is a bijective map between $\ell^2(\ZZ)$ and $\bd L^2(\TT)$, where $\TT = \RR/\ZZ$ is often identified with the interval $D_d = (-1/2,1/2]$. 

One straightforward method to construct warped filter banks on $\ell^2(\ZZ)$ would be to identify $\bd L^2(\TT)$ with $\bd L^2(D_d)$ and consider a warping function $\Phi : D \rightarrow \RR$ for some finite interval $D\subset D_d$. The elements of the warped filterbank are obtained as in Defintion \ref{def:warpedfilterbanks}, but using the inverse DTFT $\mathcal{F}_d^{-1}$ in \eqref{eq:warpedfilterbank} instead of $\mathcal{F}^{-1}$ and integer decimation factors $\bd a = (a_m)_{m\in\ZZ}\subset \NN$. When constructing warped filter banks in this fashion, frame-theoretic properties of the discrete warped filter bank are inherited from its continuous counterpart. Note that rational decimation factors can be achieved~\cite{kovacevic1993perfect}, but will not be considered here. 

In practice, the construction outlined above has a significant drawback: Despite the frequency domain being finite, the number of filters $g_m$ is still infinite. It is more common, e.g., in the construction of constant-Q \cite{brown1991calculation} and gammatone \cite{patterson1987efficient} filter banks, to restrict the desired frequency scale to an \emph{essential range} $\Xi = [\xi_{\textrm{min}},\xi_{\textrm{max}}]\subsetneq D_d$, outside of which processing, i.e., modification,  of the signal content is not desired. \revised{Then, a filter $g_m$ is only considered when its essential support has significant overlap with $\Xi$. Filters $g_m$ with essential support partially or fully outside $D_d$ are discarded. However, the filter banks so constructed are badly conditioned, or do not form frames for $\bd L^2(\TT)$ at all. More specifically the lower bound $A$ in the necessary condition 
\[
 0 < A \leq \sum_m \frac{1}{a_m}|g_m(\xi)|^2 \leq B < \infty, \text{ for all } \xi\in D_d,
\]
see Proposition \ref{pro:necandsuf}(ii), is very small, such that $A\ll B$, or $A=0$.} The above necessary condition can be satisfied by designing additional filter(s) to cover the range $D_d\setminus \Xi$, see \cite{holighaus2013framework,necciari18}.

In the following, we describe one possible method to obtain discrete warped filter banks by restricting a warping function $\Phi$ to (a subset of) $D_d$, a summary and update of work previously presented in~\cite{hoprwi15}. Although not strictly necessary, we will make the additional simplification to consider a prototype $\theta$ with $\supp(\theta)\subset [c,d]$ and $\Phi^{-1}(\Xi + [c,d]) \subset D_d$. With this assumption, we can forgo the cumbersome treatment of discontinuities and other undesired effects at the periodic boundary of $\TT$. 
Let $\Phi:D\rightarrow \RR$ be a warping function and define
\begin{equation*}
     \begin{split}
        m_{\textrm{max}} = \max \{ m\in\ZZ  \colon  \Phi^{-1}(m) \in \Xi\}+1\\
        m_{\textrm{min}} = \min \{ m\in\ZZ  \colon  \Phi^{-1}(m) \in \Xi\}-1.
     \end{split}
\end{equation*}
For $m_{\textrm{min}}< m < m_{\textrm{max}}$ the frequency responses \revised{$g_m\in \LD$ are given by the trivial restriction 
\begin{equation}
  g_m(\xi) :=  \sqrt{a_m} (\bd T_m\theta)\circ \Phi(\xi), \quad \text{ for a.e. } \xi\in D_d.
\end{equation}}
Assume that 
\begin{equation}\label{eq:fullsum}
 0 < A \leq \sum_{m\in\ZZ} |\bd T_m\theta|^2 \leq B < \infty, \text{ a.e. on } \RR.
\end{equation}
Then, it is expected that 
\begin{equation}\label{eq:restsum}
 0 < \tilde{A} \leq \sum_{m=m_{\textrm{min}}+1}^{m_{\textrm{max}}-1}  \frac{1}{a_m}|g_m|^2, \text{ a.e. on } \Xi,
\end{equation}
for some $\tilde{A}\approx A$. On $D_d\setminus \Xi$, however, the restricted sum in \eqref{eq:restsum} decays to $0$. This deficiency can be overcome if we define  
\begin{equation}
   \begin{split}
     {g_{m_{\text{min}}}}(\xi) &:= \left(a_{m_{\text{min}}-1}  \sum_{m<
   m_{\text{min}}} |(\bd T_m\theta)\circ \Phi(\xi)|^2\right)^{1/2}\quad \text{ and}\\
      {g_{m_{\text{max}}}}(\xi) &:= \left(a_{m_{\text{max}} +1}  
\sum_{m>m_{\text{max}}} |(\bd T_m\theta)\circ \Phi(\xi)|^2\right)^{1/2}, \text{ for all } \xi\in D_d,      
     \end{split}
 \label{eq:discwarpedcase2b}
\end{equation}
for some $a_m\in\NN$. Then, with the constants $A,B$ as in \eqref{eq:fullsum},
\begin{equation}\label{eq:finalsum}
 0 < A \leq \sum_{m=m_{\textrm{min}}}^{m_{\textrm{max}}}  \frac{1}{a_m}|g_m|^2 \leq B < \infty, \text{ a.e. on } D_d.
\end{equation}
As in the continuous case, the decimation factors $a_m$, and consequently the normalization of filters $g_m$, are often fixed post-hoc. \revised{Note that, inserting warped filters in \eqref{eq:finalsum}, the decimation factors cancel} and thus do not influence the bounds $A,B$. The role of $g_{m_{\text{min}}},\ g_{m_{\text{max}}}$ is to preserve the signal information in the frequency range $D_d\setminus \Xi$. The proposed definition \eqref{eq:discwarpedcase2b} prioritizes \eqref{eq:finalsum} over other properties of $g_{m_{\text{min}}},\ g_{m_{\text{max}}}$, e.g., smoothness at the periodic boundary of $\TT$. 
Clearly, it is possible to restore the lower bound in \eqref{eq:finalsum} with other filter constructions and/or a refinement of \eqref{eq:discwarpedcase2b}. Nonetheless, \eqref{eq:discwarpedcase2b} provides a straightforward method to preserve the bounds of \eqref{eq:fullsum}.
 
After the obvious changes, the necessary and sufficient conditions in Proposition \ref{pro:necandsuf} translate one-to-one to the discrete case, such that the property \eqref{eq:finalsum} is the crucial step to obtain warped filter bank frames for $\ell^2(\ZZ)$. In particular, we can apply Proposition \ref{pro:necandsuf}\revised{(iii)} to find that, whenever \eqref{eq:finalsum} holds and 
\begin{equation}\label{eq:discPL1} 
   a_m^{-1} \geq \Phi^{-1}(d+m)-\Phi^{-1}(c+m), \text{ for all } m_{\textrm{min}}< m < m_{\textrm{max}},
\end{equation}
\begin{equation}\label{eq:discPL2}
   \begin{split}
    a_{m_{\textrm{min}}}^{-1} \geq \Phi^{-1}(d+m_{\textrm{min}}) + 1/2, \text{ and} \\
    a_{m_{\textrm{max}}}^{-1} \geq 1/2-\Phi^{-1}(c+m_{\textrm{max}}), 
    \end{split}
\end{equation}
\revised{then} the discrete filter bank obtained from $(g_m)_{m\in[m_{\textrm{min}},m_{\textrm{max}}]}$ and $(a_m)_{m\in[m_{\textrm{min}},m_{\textrm{max}}]}$ is a frame with bounds $A,B$. If the decimation factors $a_m$ do not satisfy \eqref{eq:discPL1} or \eqref{eq:discPL2}, then \eqref{eq:finalsum} is no longer a sufficient frame condition. But if $\theta$ is continuous, then the frame bounds depend continuously on the choice of $a_m$. This can be seen in the discrete analogue of Proposition \ref{pro:jaleProp} (derived again from \cite[Proposition 3.7]{jale14}):
    \begin{equation}\label{eq:diagdom}    
    \infty > \tilde{B} \geq \sum_{m=m_{\textrm{min}}}^{m_{\textrm{max}}} \frac{1}{a_m} |g_m|^2 \pm \sum_{m=m_{\textrm{min}}}^{m_{\textrm{max}}}\left( 
       \frac{1}{a_m} |g_m| \sum_{l=1}^{a_m-1} |\bd T_{l/a_m} \overline{g_m}|\right) 
\geq \tilde{A} > 0, \text{ a.e. on } D_d
  \end{equation}
  is sufficient for the frame property. \revised{Note that in \eqref{eq:diagdom} (and only there), $\bd T_{x}$ denotes circular translation on $\TT$, such that the inner sum contains exactly $a_m-1$ shifted terms.} If the sampling density is to be reduced, we suggest choosing $a_m = \lfloor  \beta(\Phi^{-1}(d+m)-\Phi^{-1}(c+m))^{-1} \rfloor$, for some $\beta>1$ and $m_{\textrm{min}}< m < m_{\textrm{max}}$. This sampling scheme can be motivated, see Corollary \ref{cor:PLnatural} and Theorem \ref{thm:wpdsuffcond}, as being close to natural decimation factors. If $\Xi \subsetneq D_d$, then \revised{increasing the decimation factors $a_{m_{\textrm{min}}},a_{m_{\textrm{max}}}$ above the bounds} in \eqref{eq:discPL2} leads to quickly deteriorating frame bounds and must be handled with care.

\begin{figure}[t!]
\begin{small}
\begin{subfigure}{.48\linewidth}
\begin{centering}
\includegraphics[width=0.5\linewidth,height=0.5\linewidth]{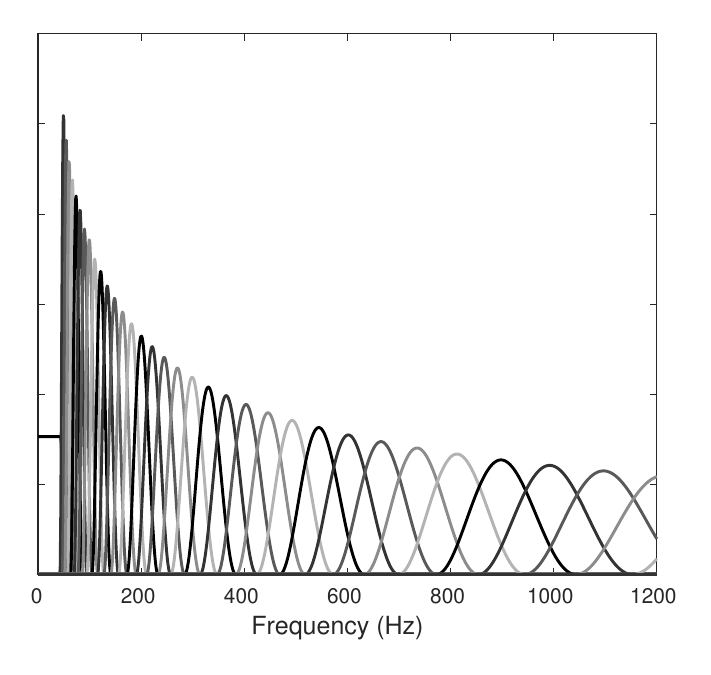} \hspace{-5pt}\includegraphics[width=0.5\linewidth,height=0.5\linewidth]{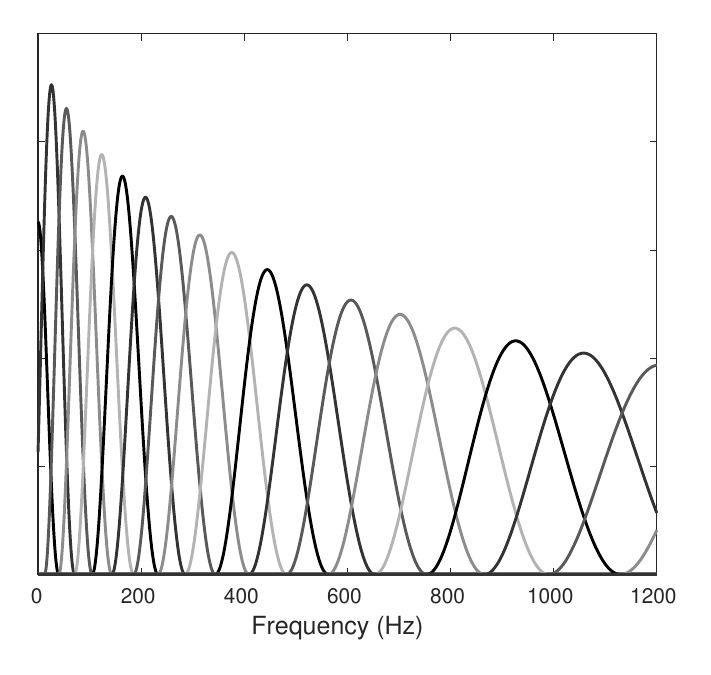}
\par\end{centering}

\begin{centering}
\includegraphics[width=0.5\linewidth,height=0.5\linewidth]{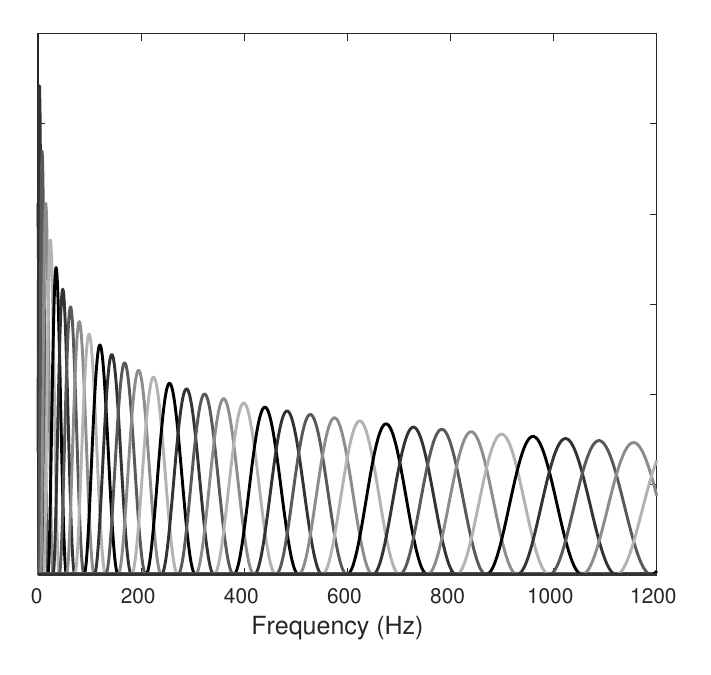} \hspace{-5pt}\includegraphics[width=0.5\linewidth,height=0.5\linewidth]{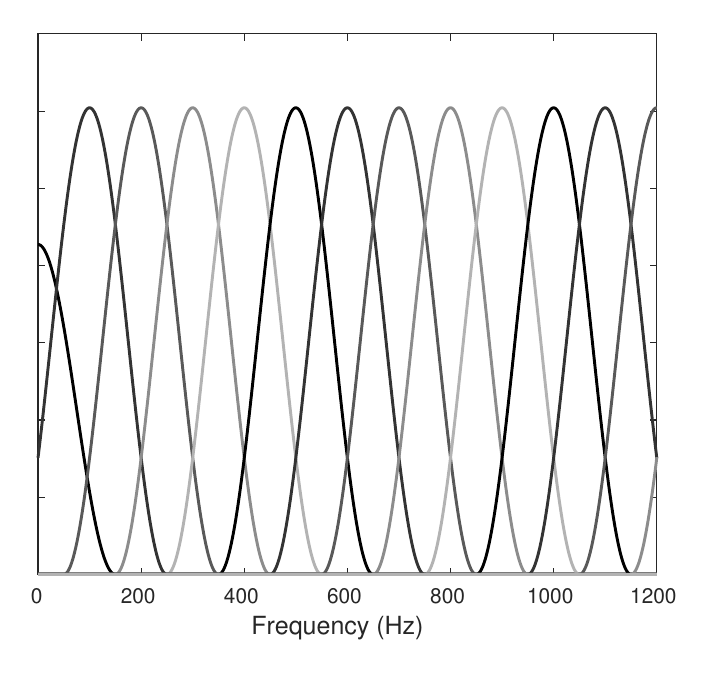}
\par\end{centering}

\caption{\label{fig:shiftinv}Frequency responses of warped filters (with low-pass filter $g_{\textrm{min}}$, see \eqref{eq:discwarpedcase2b}).
The visualization was restricted to the frequency range $0$\,Hz--$1.2$\,kHz.}
\end{subfigure}\hspace{.04\linewidth}
\begin{subfigure}{.48\linewidth} 
 \begin{centering}
\includegraphics[width=0.5\linewidth,height=0.5\linewidth]{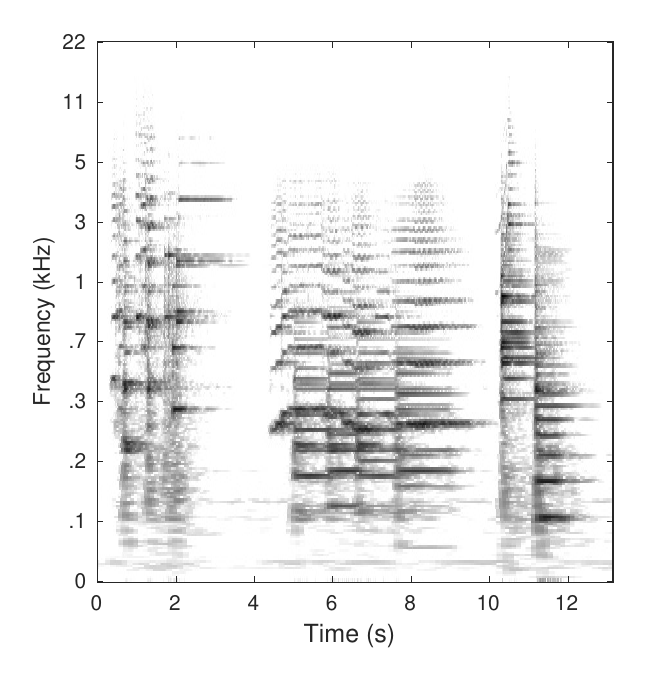} \hspace{-5pt}\includegraphics[width=0.5\linewidth,height=0.5\linewidth]{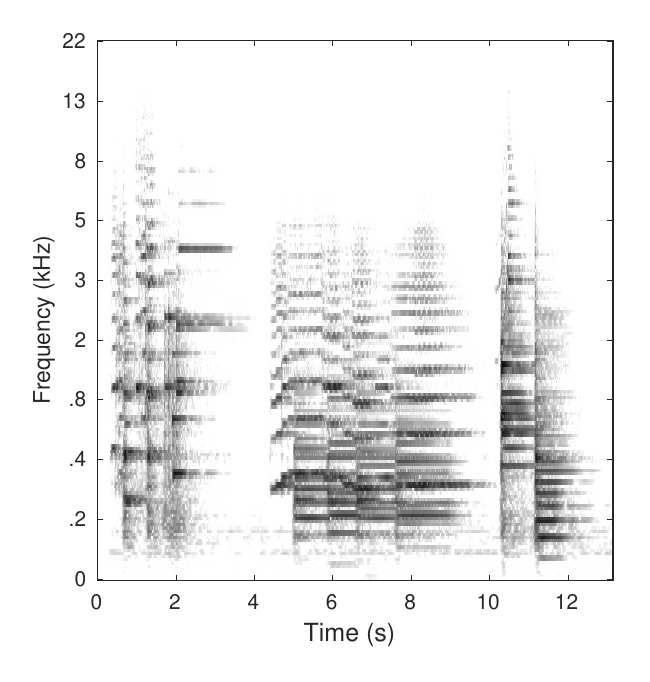}
\par\end{centering}

\begin{centering}
\includegraphics[width=0.5\linewidth,height=0.5\linewidth]{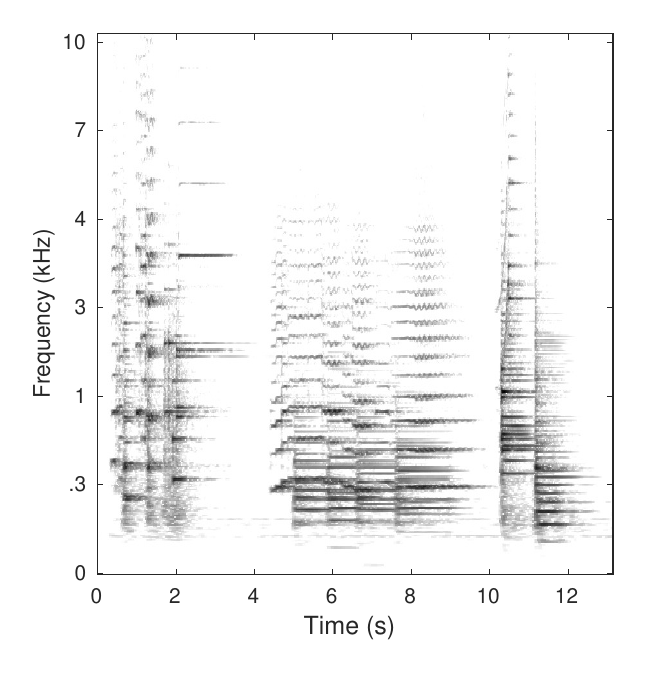} \hspace{-5pt}\includegraphics[width=0.5\linewidth,height=0.5\linewidth]{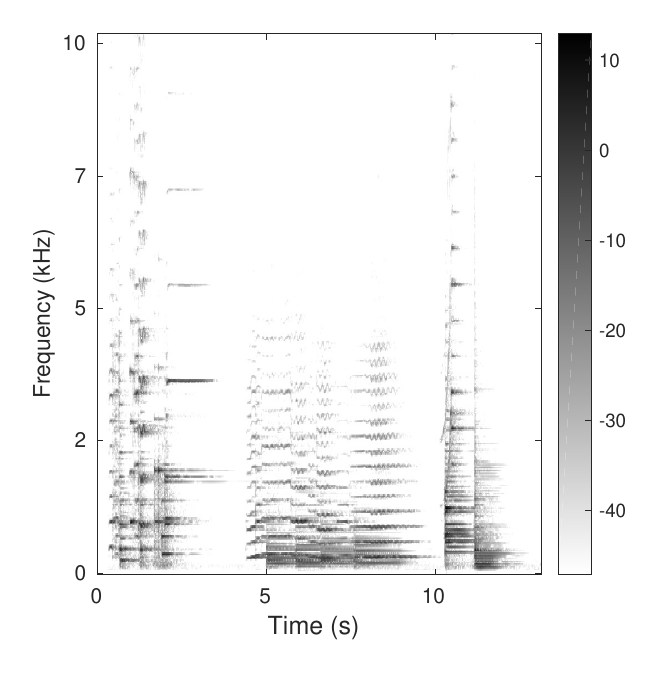}
\par\end{centering}

\caption{\label{fig:spectrograms} Time-frequency plots of a
     short piano and violin excerpt. Color indicates intensity is in dB, the colorbar is valid for all plots.}
  
\end{subfigure}
\end{small}
\caption{Warped filter bank examples for sampling rate $\xi_s = 44100$: (top-left) $\Phi_\text{log}(\xi)=10\log(\xi/\xi_s)$, (top-right) $\Phi_\text{erb}(\xi)=9.265\sgn(\xi)\log(1 + |\xi/\xi_s|/228.8)$, (bottom-left) $\Phi_\text{sqrt}(\xi)=\sgn(\xi)(\sqrt{1+|\xi/\xi_s|}-1)$, (bottom-right) $\Phi_\text{lin}(\xi)=\xi/100$. Placement applies to subfigures (a) and (b). For subfigure (b), warping functions were scaled using \eqref{eq:scadil} with $a=4$, to increase the filter density.\label{fig:2}}
\end{figure}

\revised{For some} exemplary frequency responses, derived from the warping functions introduced in Examples 
\ref{ex:wavelet1}--\ref{ex:alpha1}, see Figure \ref{fig:shiftinv}. Note that $\Phi_\text{sqrt}$ 
corresponds to Example \ref{ex:alpha1} with $\alpha=1/2$. 
In Figure \ref{fig:spectrograms}, we show time-frequency plots of a 
test signal with respect to the same warping functions.

  \subsection{Experiment: Frame bound ratio and estimates for varying redundancy}\label{ssec:ExpRatios}
  
  The results in this section, as well as Figure \ref{fig:2}, can be reproduced using the code provided at \url{http://ltfat.github.io/notes/049/}.

  The redundancy of a filter bank for $\ell^2(\ZZ)$ is given by $C_{\textrm{red}} = \sum_{m} a_m^{-1}$, where the sum is over all frequency channels. A filter bank is \emph{critically sampled} if $C_{\textrm{red}} = 1$ and
  \emph{oversampled} if $C_{\textrm{red}} > 1$\revised{, see \cite{bofehl98}. At $C_{\textrm{red}} < 1$, the filter bank is \emph{undersampled} and can never be a frame \cite{ja98,LiNg97,hoang1989non}}. However, excessive oversampling is usually undesired, as well. Warped filter bank frames may require considerable oversampling when decimation factors according to \eqref{eq:discPL1} and \eqref{eq:discPL2} are chosen. To illustrate this, assume $\supp(\theta) = [c,d]$ and $\theta(\tau) > 0$ for all $\tau\in (c,d)$. It is easy to see that \eqref{eq:discPL1} and \eqref{eq:discPL2} \revised{imply $C_{\textrm{red}} \geq \sum_m \mu(\supp(g_m))$. Furthermore, $\sum_m \mu(\supp(g_m))$ quantifies the average number of overlapping filters, since $\mu(\TT) = 1$. Finally, considering the construction of $g_m$ from integer translates of $\theta$, we obtain $\sum_m \mu(\supp(g_m))\approx \mu(\supp(\theta)) = d-c$}, i.e., the redundancy is approximately bounded below by the support size of the prototype $\theta$. From the examples in Section \ref{sec:tight}, we can see that $d-c \geq 3$ is not unexpected. 
  
  To illustrate that the redundancy can be significantly reduced below $C_{\textrm{red}} \approx \mu(\supp(\theta))$ without losing the frame property, we computed the frame bound ratio for different warped filter banks with varying redundancy, see Table \ref{tab:framebounds}. Additionally, frame bound ratio estimates in the style of  \eqref{eq:diagdom} are provided in parentheses. For every tested condition, we chose a Hann prototype, see Section \ref{sec:tight}, for $\theta$. More specifically, we chose $\vartheta$ as in  \eqref{eq:window_form} with $c_0=c_1=1/2$ and $\theta = \vartheta(\cdot/3)$. Hence, $\sum_m |\bd T_m \theta|^2 \equiv 9/8$ by \eqref{eq:sumSquares}. The first column represents filter banks with decimation factors minimizing \eqref{eq:discPL1} and \eqref{eq:discPL2} (i.e. $C_{\textrm{red}}\approx 3$). Furthermore, we tested at redundancies $2,\ 3/2,\ 5/4$ and $9/8$. Decimation factors for the redundancies below $3$ were chosen according to the sampling scheme introduced after \eqref{eq:diagdom}, with some \revised{$\beta > 1$} such that $C_{\textrm{red}} = \sum_{m} a_m^{-1}$ matched the desired redundancy closely. Even for redundancy as low as $9/8$, the frame bound ratio\footnote{In fact, due to a bug in older versions of the LTFAT Toolbox (\url{ltfat.github.io}) used for the frame bound calculations, the frame bound ratios reported in \cite{hoprwi15} are too large. In Table \ref{tab:framebounds}, corrected values are shown alongside ratios for lower redundancies not tested before.} is significantly smaller than $10$. Considering the estimate \eqref{eq:FinalPEstimate} in the proof of Lemma \ref{lem:offdiagest}, it is noteworthy that, for fixed redundancy, the dependence of the frame bound ratio on the warping function $\Phi$ seems to be quite limited.
  
  Complementing the numerically obtained ratios in Table \ref{tab:framebounds}, we computed the estimates used to prove Theorem \ref{thm:wpdsuffcond}, with $\Phi_\text{log}(\xi) = 10\log(\xi)$ and $\theta$ a Hann prototype with $R=3$ as before. Since $\theta$ is compactly supported, $\epsilon$ in the estimate \eqref{eq:FinalPEstimate} is arbitrary, but the constant $C_1C_2$ changes with $\epsilon$. We only considered the setting $a_m=\widetilde{a_w}/w(m)$ with $\widetilde{a_w}$ as in Corollary \ref{cor:PLnatural}, where we know that in fact $P_{\Phi,\theta,\bd a} \equiv 0$, see \eqref{eq:defOfP} for the definition of $P_{\Phi,\theta,\bd a}$. For $\epsilon = 1,2,3,4,5$, Equation \eqref{eq:FinalPEstimate} yields for $P_{\Phi,\theta,\bd a}$ the upper bounds $10.2,~6.5,~6.6,~8.2,~11.6$ (rounded down to the first decimal), such that Theorem \ref{thm:wpdsuffcond} would 
not be sufficient to confirm the frame property, although the considered system is even tight. It might be interesting to note that, for small $\epsilon$, the dominant quantity in \eqref{eq:FinalPEstimate} is the first term in parentheses, while for larger $\epsilon$, the constant $C_1C_2$ is dominant. Curiously, the term depending on $A_v^{-1}$ had relatively minor contribution.
%
%

\begin{table}[t!] 
    \caption{Frame bound ratios of warped filter banks from Figure~\ref{fig:shiftinv} with 
    varying redundancy. Columns correspond to warped filter banks with \emph{approximately} equal
    redundancy. Numbers in parentheses are estimates obtained by considering the sum and difference, respectively, of the terms in \eqref{eq:diagdom}.}\label{tab:framebounds}

\begin{center}\begin{tabular}{ |c|ccccc|} 
 \hline   
     $C_{\text{red}}(\approx)$ &         $3$ &                  $2$ &               $3/2$ &                $5/4$ &           $9/8$ \\ 
 \hline 
     $\Phi_\text{lin}$         & $   1.000 $ & $   1.220 $ ($1.234$)& $   1.961 $ ($1.982$)& $   3.880 $ ($4.759$)& $6.868$ ($10.042$)\\ 
     $\Phi_\text{sqrt}$        & $   1.003 $ & $   1.237 $ ($1.243$)& $   1.980 $ ($1.997$)& $   3.938 $ ($4.894$)& $7.315$ ($11.135$)\\ 
     $\Phi_\text{erb}$         & $   1.000 $ & $   1.240 $ ($1.249$)& $   1.970 $ ($2.134$)& $   3.860 $ ($5.023$)& $7.122$ ($11.323$)\\ 
     $\Phi_\text{log}$         & $   1.014 $ & $   1.240 $ ($1.249$)& $   1.973 $ ($2.125$)& $   3.876 $ ($5.019$)& $7.159$ ($11.323$)\\ 
 \hline 
 \end{tabular}\end{center} 
\end{table}   
%

\section{Conclusion and Outlook}
  
  In this contribution, we have introduced a novel, flexible family of structured time-frequency filter banks. These warped filter banks are able to recreate (or imitate) important classical time-frequency representations, while providing additional design freedom. Warped filter banks allow for intuitive handling and the application of important results from the theory of generalized-shift invariant frames. In particular, the construction of tight frames of bandlimited filters is easy: It reduces to the selection of a compactly supported prototype function whose integer translates satisfy a simple summation condition and sufficiently small decimation factors $a_m$. Moreover, the warping construction induces a natural choice of decimation factors that further simplifies the design of warped filter bank frames. With several examples, we have illustrated not only the flexibility of our method when selecting a non-linear frequency scale, but also the ease with which tight or snug frames can be constructed.

  The complementary manuscript~\cite{bahowi14} discusses warped time-frequency representations in the context of continuous frames, determines the associated coorbit spaces and the warped time-frequency representations' sampling properties in the context of atomic decompositions and Banach frames~\cite{fegr89,fegr89-1}. Future work will continue to explore practical applications of warped time-frequency representations and their discrete equivalents, as well as extending the warping method to multidimensional signals.
  
\section*{Acknowledgment}
  This work was supported by the Austrian Science Fund (FWF)
  projects FLAME (Y 551-N13) and MERLIN (I 3067-N30) and the Vienna Science and Technology Fund (WWTF) Young Investigators project CHARMED (VRG12-009). The authors wish to thank Felix Voigtlaender and Peter Balazs for fruitful discussion on the topics covered here, Jakob Lemvig and G\"unther Koliander for important constructive comments on preprints of the manuscript and the anonymous reviewers for their tremendously helpful comments and suggestions.


 \bibliographystyle{abbrv}
 \section*{\refname}
 

 \end{document}